\documentclass{amsart}
\title[CMC-1 surfaces in de Sitter space]{
    Spacelike mean curvature one surfaces \\
    in de Sitter $3$-space%
}
\usepackage[dvips]{graphicx}
\usepackage{enumerate}
\usepackage{amssymb}
\usepackage{amsthm}
\theoremstyle{plain}
 \newtheorem{theorem}{Theorem}[section]
 \newtheorem*{theorem*}{Theorem}
 \newtheorem{introtheorem}{Theorem}
   
 \newtheorem*{lemma*}{Lemma}
 \newtheorem*{lemmaE1}{Lemma~E1}
 \newtheorem*{lemmaE2}{Lemma~E2}
 \newtheorem*{lemmaP}{Lemma~P}
 \newtheorem*{lemmaH}{Lemma~H}
 \newtheorem{proposition}[theorem]{Proposition}
 \newtheorem{fact}[theorem]{Fact}
 \newtheorem*{fact*}{Fact}

 \newtheorem{lemma}[theorem]{Lemma}
 \newtheorem{corollary}[theorem]{Corollary}
 \theoremstyle{remark}
 \newtheorem{definition}[theorem]{Definition}
 \newtheorem{remark}[theorem]{Remark}
 \newtheorem*{remark*}{Remark}
 \newtheorem*{problem*}{Problem}
 
 \newtheorem{example}[theorem]{Example}
\numberwithin{equation}{section}
\numberwithin{figure}{section}

\newcommand{\Z}{\boldsymbol{Z}}
\newcommand{\R}{\boldsymbol{R}}

\newcommand{\C}{\boldsymbol{C}}
\newcommand{\Herm}{\operatorname{Herm}}

\newcommand{\SL}{\operatorname{SL}}
\newcommand{\PSL}{\operatorname{PSL}}
\newcommand{\SU}{\operatorname{SU}}
\newcommand{\SO}{\operatorname{SO}}
\renewcommand{\O}{\operatorname{O}}
\newcommand{\PSU}{\operatorname{PSU}}

\newcommand{\LC}{\operatorname{\mathit{LC}}}
\newcommand{\trace}{\operatorname{trace}}
\newcommand{\Ord}{\operatornamewithlimits{Ord}}
\newcommand{\id}{\operatorname{id}}
\newcommand{\sign}{\operatorname{sgn}}
\newcommand{\inner}[2]{\left\langle{#1},{#2}\right\rangle}

\renewcommand{\Re}{\operatorname{Re}}
\renewcommand{\Im}{\operatorname{Im}}

\renewcommand{\phi}{\varphi}
\renewcommand{\epsilon}{\varepsilon}

\newcommand{\cmcone}{\mbox{CMC-$1$}}
\newcommand{\secondff}{\mbox{I\!I}}
\def\transpose#1{\mathord{\mathopen{{\vphantom{#1}}^t}#1}}

\DeclareFontFamily{U}{rsfs}{\skewchar\font"7F}
\DeclareFontShape{U}{rsfs}{m}{n}{
        <-6> rsfs5
        <6-8> rsfs7
        <8-> rsfs10
        }{}
\DeclareMathAlphabet{\mathscr}{U}{rsfs}{m}{n}
\author{S.~Fujimori}
\address[Shoichi Fujimori]{%
   Department of Mathematics, Fukuoka University of Education,
   Munakata, Fukuoka 811-4192, Japan}
\email{fujimori@fukuoka-edu.ac.jp}

\author{W.~Rossman}
\address[Wayne Rossman]{%
   Department of Mathematics, Faculty of Science,
   Kobe University,
   Rokko, Kobe 657-8501, Japan
}
\email{wayne@math.kobe-u.ac.jp}

\author{M.~Umehara}
\address[Masaaki Umehara]{%
   Department of Mathematics, Graduate School of Science,
   Osaka University,
   Toyonaka, Osaka 560-0043,
   Japan
}
\email{umehara@math.sci.osaka-u.ac.jp}

\author{K.~Yamada}
\address[Kotaro Yamada]{%
   Faculty of Mathematics,
   Kyushu University,
   Higashi-ku, Fukuoka 812-8581, Japan
}
\email{kotaro@math.kyushu-u.ac.jp}

\author{S.-D.~Yang}
\address[Seong-Deog Yang]{%
   Department of Mathematics,
   Korea University,
   Seoul 136-701, Korea
}
\email{sdyang@korea.ac.kr}


\begin{document}
\begin{abstract}
 The first author studied
 spacelike constant mean curvature one (\cmcone)
 surfaces in de Sitter 3-space $S^3_1$
 when the surfaces have no singularities except within some
 compact subset and are of finite total curvature on
 the complement of this compact subset.
 However, there are many \cmcone{} surfaces whose
 singular sets are not compact.
 In fact, such examples have already appeared in the
 construction of trinoids given by Lee and the last author
 via hypergeometric functions.

 In this paper, we improve the Osserman-type inequality
given by the first author.
Moreover, we shall develop a fundamental framework that allows
 the singular set to be non-compact, and then
 will use it to investigate the global behavior of \cmcone{} surfaces.
\end{abstract}
\maketitle

\section*{Introduction}
A holomorphic map $F:M^2\to \SL_2\C$ of  a Riemann surface
$M^2$ into the complex Lie group $\SL_2\C$ is called {\em null\/}
if $\det(dF/dz)$ vanishes identically, where $z$ is a local complex
coordinate of $M^2$.
We consider two projections, one into the hyperbolic $3$-space
\[
  \pi_H:\SL_2\C\longrightarrow H^3=\SL_2\C/\SU_2
\]
and the other into the de Sitter $3$-space
\[
  \pi_S:\SL_2\C\longrightarrow S^3_1=\SL_2\C/\SU_{1,1},
\]
where the definition of $\SU_{1,1}$ is in Appendix~\ref{app:conj}.
It is well-known that the projection of a holomorphic null immersion
into $H^3$ by $\pi_H$ gives a conformal \cmcone{}
(constant mean curvature one) immersion
(see \cite{Br}, \cite{UY0}, \cite{CHR}).
Moreover, conformal \cmcone{} immersions are always given locally
in such a manner.

On the other hand, spacelike \cmcone{} surfaces given by
the projection of holomorphic null immersions
into $S^3_1$ by $\pi_S$ can have singularities,
and are called {\em \cmcone{} faces}. 
We work with this class of surfaces that is larger than the class of
\cmcone{} immersions. 
In fact, the class of \cmcone{} immersions is too 
small,
since {\it there is only one, up to congruency, 
complete spacelike \cmcone{}
immersion} \cite{Ak,R}, which we call an {\em $S^3_1$-horosphere}.
(We also give a simple proof of this here. See the last remark of 
Section~\ref{sec:prelim}.)

The relationship between \cmcone{} surfaces in $H^3$ and
\cmcone{} faces in $S^3_1$ is analogous to that between
minimal surfaces in Euclidean 3-space $\R^3$ and spacelike 
maximal surfaces with singularities in Lorentz-Minkowski 3-space
$\R_1^3$ (called {\em maxfaces} \cite{UY1}). 
Note that maximal surfaces also admit a Weierstrass-type representation
formula (\cite{K}).
As in the case of maxfaces (see \cite{UY1}),
the first author \cite{F} investigated the global behavior of
\cmcone{} faces in $S^3_1$, in particular
proving an Osserman-type inequality for complete
\cmcone{} faces of finite type whose ends are all elliptic,  
where a complete end of a \cmcone{} face is called
{\it elliptic, parabolic}, or {\em hyperbolic\/} if the monodromy matrix
of the holomorphic lift $F:M^2\to \SL_2\C$ is elliptic, parabolic, 
or hyperbolic, respectively (see Section~\ref{sec:prelim}).  
One of our main results is the following, which implies that 
the ellipticity or parabolicity 
of ends follows from  completeness:
\begin{introtheorem}\label{thm:main}
 A complete end of a \cmcone{} face in $S_1^3$ is never hyperbolic,
 so must be either elliptic or parabolic.
 Moreover, the total curvature over a neighborhood of such an end is
 finite.
\end{introtheorem}
We remark that there exist incomplete elliptic and parabolic ends.

It is remarkable that just completeness of an end
is sufficient to conclude that it has finite total curvature.  This
is certainly not the case for \cmcone{} surfaces in $H^3$ nor for
minimal surfaces in $\R^3$, 
but is similar to the case of maximal surfaces in $\R^3_1$ \cite{UY1}.
Although the asymptotic behavior of regular elliptic \cmcone{}
ends in $S_1^3$ is investigated in \cite{F}, there do also exist
complete parabolic ends, and to describe them, a much deeper
analysis is needed, which we will conduct in this article.

As an application of Theorem \ref{thm:main}, we prove the following
Osserman-type inequality, which improves the result of \cite{F}
by removing the assumptions of finite type and ellipticity of ends:
\begin{introtheorem}
\label{thm:fujimori}
 Suppose a  \cmcone{} face $f\colon{}M^2 \to S^3_1$ is complete.
 Then there exist a compact Riemann surface $\overline M^2$ and a finite
 number of points
 $p_1,\dots,p_n \in \overline M^2$ such that $M^2$ is biholomorphic to
 $\overline M^2\setminus\{p_1,\dots,p_n\}$, and
 \begin{equation}\tag{$*$}\label{eq:fujimori}
   2\deg(G)\ge -\chi(\overline M^2)+2n, 
 \end{equation}
 where $G$ is the hyperbolic Gauss map of $f$ and $\chi(\overline M^2)$ is the
 Euler characteristic of $\overline M^2$.  Furthermore, equality
 holds if and only if
 each end is regular and properly embedded.
\end{introtheorem}

\cmcone{} trinoids in $S^3_1$ were constructed
by Lee and the last author using hypergeometric functions \cite{LY}, and
those trinoids with elliptic ends are complete in the sense of \cite{F},
and attain equality in \eqref{eq:fujimori}.
However, those having other types of ends are not complete, as
their singular sets are not compact.
For this reason, our goal is not only to prove the above two theorems,
but also to extend the framework for
\cmcone{} surfaces to include a larger class of surfaces,
relaxing the immersedness and completeness conditions.
If $M$ is of finite topology, i.e. if $M$ is diffeomorphic to a compact
Riemann surface $\overline M^2$ with finitely many punctures 
$p_1,\dots, p_n$,
and if a \cmcone{} face $f:M^2 \to S^3_1$ is
weakly complete, whose precise definition will  be given in
Section~\ref{sec:prelim},
we say that $f$ is a weakly complete \cmcone{} face of finite topology.
We shall develop the framework under this more general notion, which includes
all the trinoids in \cite{LY}.

In Section~\ref{sec:prelim}, we recall definitions and basic
results.
In Section~\ref{sec:projective}, we investigate the
monodromy of the hyperbolic metrics
on a punctured disk around an end.
As an application, we prove Theorem \ref{thm:main}
in Section~\ref{sec:ends}.
In Section~\ref{sec:lgm}, we give a geometric interpretation of the
hyperbolic Gauss map.
In Section~\ref{sec:fujimori}, we investigate the asymptotic
behavior of regular parabolic ends, and prove Theorem
\ref{thm:fujimori}.
In Appendix~\ref{app:Q-mero}, we prove meromorphicity of the Hopf
differential for complete \cmcone{} faces.
In Appendix~\ref{app:conj}, we explain the conjugacy classes of $\SU_{1,1}$.

Generic singularities of \cmcone{} faces are classified in \cite{FSUY}.
A \cmcone{} face is called {\it embedded\/}
(in the wider sense) if it is embedded outside of some 
compact set of $S^3_1$.  Examples of complete embedded 
\cmcone{} faces are given in \cite{FRUYY} as deformations 
of maxfaces in Lorentz-Minkowski 3-space
$\R_1^3$.

\section{Preliminaries}
\label{sec:prelim}
\subsection*{The representation formula}
Let $\R^4_1$ be the Lorentz-Minkowski
space of dimension $4$, with the Lorentz metric
\[
  \langle (x_0,x_1,x_2,x_3),(y_0,y_1,y_2,y_3)\rangle =
  -x_0y_0+x_1y_1+x_2y_2+x_3y_3.
\]
Then de Sitter $3$-space is
\[
  S^3_1=
  \{(x_0,x_1,x_2,x_3)\in\R^4_1 \, ; \, -x_0^2+x_1^2+x_2^2+x_3^2=1\},
\]
with metric induced from $\R^4_1$, which is
a simply-connected Lorentzian $3$-manifold with constant
sectional curvature $1$.
We identify $\R^4_1$ with the set of 
$2\times 2$ Hermitian matrices
$\Herm(2)=\{X^* = X\}$ ($X^*:=\transpose{\overline{X}}$) by
\begin{equation}\label{eq:herm-mink}
   X=(x_0,x_1,x_2,x_3)\leftrightarrow
    X=\sum_{k=0}^3 x_k e_k
    =\begin{pmatrix}
      x_0+x_3           & x_1+i x_2 \\
      x_1-i x_2 & x_0-x_3
     \end{pmatrix},
\end{equation}
where 
\begin{equation}\label{eq:pauli}
    e_0=\begin{pmatrix}1&0\\0&1\end{pmatrix},~
    e_1=\begin{pmatrix}0&1\\1&0\end{pmatrix},~
    e_2=\begin{pmatrix}\hphantom{-}0&i\\-i&0\end{pmatrix},~
    e_3=\begin{pmatrix}1&\hphantom{-}0\\0&-1\end{pmatrix}
\end{equation}
and $i=\sqrt{-1}$. 
Then $S^3_1$ is
\[
  S^3_1=\{X\,;\,X^*=X\,,\det X=-1\}
              =\{Fe_3F^*\,;\,F\in \SL_2\C\}
\]
with the metric
\[
   \inner{X}{Y}
    =-\frac{1}{2}\trace\left(Xe_2(\transpose{Y})e_2\right),
    \qquad
   \inner{X}{X}=-\det X.
\]
The projection $\pi_S\colon{}\SL_2\C\to S^3_1$ mentioned in the introduction 
is written explicitly as $\pi_S(F)=F e_3 F^*$.
Note that the hyperbolic $3$-space $H^3$ is given by 
$H^3=\{FF^*\,;\,F\in\SL_2\C\}$ and the projection is
$\pi_H(F)=FF^*$.

An immersion into ${S}^3_1$ is called {\em spacelike\/}
if the induced metric on the immersed surface is positive definite.
The complex Lie group $\SL_2\C$ acts isometrically on
$\Herm(2)=\R^4_1$, as well as $S^3_1$, by
\begin{equation}\label{eq:isometry-action}
   \Herm(2)\ni X \longmapsto a X a^*\qquad a\in\SL_2\C.
\end{equation}
In fact, $\PSL_2\C=\SL_2\C/\{\pm\id\}$ is isomorphic
to the identity component $\SO^+_{3,1}$ of the isometry group
$\O_{3,1}$ of $S^3_1$.
Note that each element of $\SO^+_{3,1}$ corresponds to an
orientation preserving and orthochronous 
(i.e., time orientation preserving)
isometry.
The group of orientation preserving isometries of $S^3_1$
is generated by $\PSL_2\C$ and
the map
\begin{equation}\label{eq:time-reverse}
   S^3_1\ni X \longmapsto -X \in S^3_1.
\end{equation}

Aiyama-Akutagawa \cite{AA} gave a Weierstrass-type representation
formula in terms of holomorphic data 
for spacelike \cmcone{} immersions in $S^3_1$. 
The first author \cite{F} extended 
the notion of \cmcone{} surfaces as follows, 
  like as for the case of maximal surfaces in the Minkowski space \cite{UY1}.
\begin{definition}[\cite{F}]
 Let $M^2$ be a $2$-manifold. 
 A $C^\infty$-map $f\colon{}M^2\to {S}^3_1$ is called a 
 {\em \cmcone{} face\/} if 
 \begin{enumerate}
  \item there exists an open dense subset $W\subset M^2$ such that 
        $f|_W$ is a spacelike \cmcone{} immersion, 
  \item for any singular point 
        (that is, a point where the induced metric degenerates) $p$, 
        there exists a $C^1$-differentiable function 
        $\lambda :U\cap W\to (0,\infty)$, 
        defined on the intersection of 
        neighborhood $U$ of $p$ with $W$,        
        such that $\lambda\,ds^2$ extends to a $C^1$-differentiable 
        Riemannian metric on $U$, 
        where $ds^2$ is the first fundamental form, i.e.,
        the pull-back of the metric of $S^3_1$ by $f$,         
        and 
  \item $df(p)\ne 0$ for any $p\in M^2$. 
 \end{enumerate}
\end{definition}

\begin{remark}\label{rmk:ori}
Though the original definition of \cmcone{} faces 
in \cite{F} assumed the orientability of 
the source manifold, our definition here does not.
However, this difference is not of an essential 
nature.  In fact, 
for any \cmcone{} face $f:M^2\to S^3_1$,
$M^2$ is automatically orientable.
(See \cite{KU}.)
\end{remark}

\begin{remark}\label{rmk:causal}
A $C^\infty$-map $f:M^2\to S^3_1$ is called
a  {\it frontal} 
if $f$ lifts to a 
$C^\infty$-map $L_f \colon M^2 \to P(T^*S^3_1)$
such that $dL_f(TM^2)$ lies in the
canonical contact plane-field on $P(T^*S^3_1)$.
Moreover, $f$ is called a
{\it wave front} or a {\it front} if 
$L_f$ is an immersion, that is, $L_f(M^2)$ is  
a Legendrian submanifold.
If a frontal $L_f$ can lift up to a
smooth map into $T^*S^3_1$, $f$ is called 
{\it co-orientable}, and otherwise it is called
{\it non-co-orientable}.
Wave fronts are a canonical class for 
investigating flat surfaces in the hyperbolic 3-space $H^3$.
In fact, like for \cmcone{} faces (see Theorem \ref{thm:fujimori}
in the introduction), an Osserman-type
inequality holds for flat fronts in $H^3$ (see \cite{KUY}.)    
Although our \cmcone{} faces belong to a special class of 
horospherical linear Weingarten surfaces (cf. \cite{KU}),
they may not be  (wave) fronts in general, but are
co-orientable frontals.
In particular, there is a globally defined non-vanishing normal vector field
$\nu$ on the whole of $M^2$ for a given \cmcone{} face $f:M^2\to S^3_1$.
It should be remarked that the limiting tangent plane
at each singular point contains a lightlike direction, 
that is, a \cmcone{} face is not spacelike on 
the singular set.
\end{remark}

An oriented $2$-manifold $M^2$ on which a \cmcone{} face
$f:M^2\to {S}^3_1$ is defined always has a complex structure 
(see \cite{F}). 
Since \cmcone{} faces are all
orientable and co-orientable (cf. \cite{KU}), 
from now on, we will treat 
$M^2$ as a Riemann surface, and we can assume the existence 
of a globally defined non-vanishing 
normal vector field. 
The representation formula in \cite{AA} can be extended 
for \cmcone{} faces as follows:

\begin{theorem}[{\cite[Theorem 1.9]{F}}]\label{th:AA-rep}
 Let $\widetilde{M}^2$ be a simply connected Riemann surface.
 Let $g$ be a
 meromorphic function and $\omega$ a holomorphic $1$-form on 
 $\widetilde{M}^2$ such that
 \begin{equation}\label{eq:dshat^2}
  d\hat s^2=(1+|g|^2)^2|\omega|^2
 \end{equation}
 is a Riemannian metric on $\widetilde{M}^2$ and
 $|g|$ is not identically $1$.
 Take a holomorphic immersion $F=(F_{jk}):\widetilde M^2\to \SL_2\C$ satisfying
 \begin{equation}\label{eq:ode}
    F^{-1}dF=\begin{pmatrix}g&-g^2\\1&-g\hphantom{^2}\end{pmatrix}\omega.
 \end{equation}
 Then $f:\widetilde M^2\to {S}^3_1$ defined by
 \begin{equation}\label{eq:projF}
    f=\pi_S\circ F := Fe_3F^*
 \end{equation}
 is a \cmcone{} face which is conformal away from its singularities.
 The induced
 metric $ds^2$ on $\widetilde M^2$, the second fundamental form 
 $\secondff$, 
 and the Hopf differential $Q$ of $f$ are given as follows{\rm :}
 \begin{equation}\label{eq:forms}
  ds^2 =(1-|g|^2)^2 |\omega|^2 ,\quad
  \secondff = Q + \overline{Q}+ds^2,\quad Q=\omega\,dg.
 \end{equation}
 The singularities of the \cmcone{} face occur at points where $|g|=1$.

 Conversely, for any \cmcone{} face $f:\widetilde M^2\to {S}^3_1$,
 there exist a meromorphic function $g$
 {\rm (}with $|g|$ not identically $1${\rm)} and 
 a holomorphic $1$-form $\omega$ on
 $\widetilde{M}^2$ so that $d\hat{s}^2$ is a Riemannian metric on
 $\widetilde{M}^2$ and \eqref{eq:projF} holds,
 where $F:\widetilde{M}^2\to \SL_2\C$ is an immersion satisfying
 \eqref{eq:ode}.
\end{theorem}

\begin{remark}\label{rmk:line}
By definition,  \cmcone{} faces have
dense regular sets. However, the projection of
null holomorphic immersions might not have  
dense regular sets, in general.
Such an example has been given in \cite[Remark 1.8]{F}. 
Fortunately, we can explicitly classify such 
degenerate examples, as follows: 
Let $M^2$ be a connected Riemann surface and 
$F\colon{}M^2\to \SL(2,\C)$ be a null {\em immersion}.
We assume that the set of singular points of the corresponding map
\[
     f=F{e_3}F^*\colon{}M^2\longrightarrow S^3_1
\]
has an interior point.
Then the secondary Gauss map $g$ is constant on $M^2$ and $|g|=1$.
Without loss of generality, we may assume
$g =1$.
Since $F$ is an immersion, $(1+|g|^2)^2|\omega|^2$ is positive definite.
Then $\omega\neq 0$ everywhere.
Hence for each $p\in M^2$, one can take a complex coordinate $z$ 
such that $\omega=dz$.
Then $F$ is a solution of
\[
    F^{-1}dF = \begin{pmatrix}
                  1 & -1 \\
                  1 & -1
               \end{pmatrix}dz.
\]
Without loss of generality, we may assume that
$F(0)=
\begin{pmatrix}
{1}/{2} & {1}/{2} \\
        -1 & 1
\end{pmatrix}$.
Then we have 
\[
     F=\begin{pmatrix}
        z+{1}/{2} & -z+{1}/{2} \\
        -1 & 1
        \end{pmatrix},
\]
and the corresponding map $f$ is computed as
\[
    f=Fe_3F^*=\begin{pmatrix}
               2 \Re z & -1 \\
               -1 & 0
               \end{pmatrix},
\]
whose image is a lightlike line in $S^3_1$.  Thus, we 
have shown that the image of any degenerate \cmcone{} surface is 
a part of a lightlike line.
\end{remark}

\begin{remark}\label{lem:bryant}
 Theorem~\ref{th:AA-rep} is an analogue of the Bryant representation for
 \cmcone{} surfaces in $H^3$, which explains why \cmcone{} surfaces
 in both $H^3$ and $S^3_1$ are characterized by the projections
 $\pi_H \circ F$ and $\pi_S\circ F$. 
 The \cmcone{} surfaces in $H^3$ and $S^3_1$ are both
 typical examples in the class of linear Weingarten surfaces.
 A Bryant-type representation formula for
 linear Weingarten surfaces was recently given by 
 J.~G\'alvez, A.~Mart\'\i{}nez and F.~Mil\'an \cite{GMM}.
\end{remark}
\begin{remark}\label{lem:data}
 Following the terminology of \cite{UY0}, $g$ is called 
 a {\em secondary\/} Gauss map of $f$.
 The pair $(g,\omega)$ is called {\em Weierstrass data} of $f$,
 and $F$ is called a {\em holomorphic null lift\/} of $f$.

 The holomorphic $2$-differential $Q$ as in \eqref{eq:forms} is
 called the {\em Hopf differential} of $f$.
 In analogy with the theory of \cmcone{} surfaces in $H^3$,
 the meromorphic function
\begin{equation}\label{eq:h-gauss}
           G := \frac{dF_{11}}{dF_{21}}= \frac{dF_{12}}{dF_{22}}
\end{equation}
 is called the {\em hyperbolic Gauss map}.
 A geometric meaning for the hyperbolic Gauss map is
 given in Section~\ref{sec:lgm}.
\end{remark}
\begin{remark}\label{rem:minimal-maximal}
 Corresponding to Theorem~\ref{th:AA-rep},
 a Weierstrass-type representation formula is known
 for spacelike maximal surfaces in $\R^3_1$ (\cite{K}).
 In fact, the Weierstrass data $(g,\omega)$ as in
 Theorem~\ref{th:AA-rep}
 defines null curves in $\C^3$ by
 \[
  F_0(z):=\int_{z_0}^z \bigl(-2g, 1+g^2,i(1-g^2)\bigr)\omega.
 \]
 Any maxface (see \cite{UY1} for the definition) is locally obtained as
 the real part of some $F_0$. 
 Moreover, their first fundamental forms and Hopf differentials
 are given by \eqref{eq:forms}.
 The meromorphic function $g$ can be identified with
 the Lorentzian Gauss map.
 In this case, we call the pair $(g,\omega)$
 the {\em Weierstrass data\/} of the maxface.
\end{remark}
\begin{remark}\label{rem:small}
 Let $G,g$ be meromorphic functions on a Riemann surface.  Set
 \begin{equation}\label{eq:small}
  F=\begin{pmatrix}
     G\dfrac{da}{dG}-a & G\dfrac{db}{dG}-b \\[6pt]
     \dfrac{da}{dG} & \dfrac{db}{dG}
    \end{pmatrix},\qquad
    a=\sqrt{\frac{dG}{dg}},\quad
    b=-ga.
 \end{equation}
 Then $F$ is a meromorphic null map with hyperbolic and secondary Gauss
 maps $G$ and $g$.
 Formula \eqref{eq:small} is called Small's formula
 (\cite{KUY0}, \cite{Small}).
\end{remark}
\begin{remark}\label{rem:ambiguity}
 The holomorphic null lift $F$ of a \cmcone{} face $f$
 is unique up to right-multiplication by matrices in $\SU_{1,1}$,
 that is, for each $A\in\SU_{1,1}$, the projection of $FA^{-1}$
 is also $f$.
 Under the transformation $F \mapsto FA^{-1}$, 
 the secondary Gauss map $g$ changes by
 a M\"obius transformation:
 \begin{equation}\label{eq:second-moebius}
     g \longmapsto A \star g:=\frac{A_{11}g+ A_{12}}{A_{21}g + A_{22}},
      \qquad
        A=\begin{pmatrix}
           A_{11}& A_{12} \\A_{21} & A_{22}
          \end{pmatrix}.
 \end{equation}
 The conditions $|g|=1$, $|g|>1$, $|g|<1$ are
 invariant under this transformation.

 In particular, let $f\colon{}M^2\to S^3_1$ be a \cmcone{} face
 of a (not necessarily simply connected) Riemann surface $M^2$.
 Then the holomorphic null lift $F$ is defined only on the universal
 cover $\widetilde{M}^2$ of $M^2$.
 Take a deck transformation $\tau\in\pi_1(M^2)$
 in $\widetilde M^2$.
 Since $\pi_S\circ F=\pi_S\circ F\circ \tau$,
 there exists a $\tilde\rho(\tau)\in\SU_{1,1}$
 such that
 \begin{equation}\label{eq:mono-repr}
    F\circ \tau = F\tilde\rho(\tau).
 \end{equation}
 The representation $\tilde\rho\colon{}\pi_1(M^2)\to \SU_{1,1}$
 is called the {\em monodromy representation},
 which induces a 
 $\PSU_{1,1}$-representation
 $\rho\colon{}\pi_1(M^2)\to\PSU_{1,1}=\SU_{1,1}/\{\pm 1\}$ satisfying 
 \begin{equation}\label{eq:g-repr}
       g\circ\tau^{-1} = \rho(\tau)\star g.
 \end{equation}
\end{remark}
\begin{remark}\label{rem:action-null}
 The action $F\mapsto BF$, $B \in \SL_2\C$,
 induces a rigid motion $f\mapsto B f B^*$ in $S^3_1$,
 and the isometric motion $f \mapsto -f$ as in
 \eqref{eq:time-reverse} corresponds to
 \begin{equation}\label{eq:time-reverse-lift}
     F\longmapsto F^{\natural}=F\begin{pmatrix} 0 & i \\ i & 0 \end{pmatrix}
            = iFe_1.
 \end{equation}
 The secondary Gauss map of $F^{\natural}$ is $1/g$.
 \end{remark}
\begin{remark}\label{rem:gauss}
 Let $K_{ds^2}$ be the Gaussian curvature of $ds^2$ on
 the set of regular points of $f$.
 Then
 \begin{equation}\label{eq:hyp-metric}
  d\sigma^2:= K_{ds^2}\,ds^2 =
   \frac{4\,|dg|^2}{(1-|g|^2)^2}
 \end{equation}
 is a pseudometric of constant curvature $-1$,
 which degenerates at isolated umbilic points.  We have
 \begin{equation}\label{eq:gauss}
  d\sigma^2\cdot ds^2 = 4|Q|^2.
 \end{equation}
\end{remark}
\begin{remark}\label{rem:lift-metric}
 The metric
  \begin{equation}\label{eq:lift-metric}
   ds_\#^2:=
    \bigl(1+|G|^2\bigr)^2\left| \frac{Q}{dG} \right|^2
  \end{equation}
 is induced  from the canonical Hermitian metric
 of $\SL_2\C$ via  $F^{-1}\colon{}\widetilde M^2\to \SL_2\C$.
 When the \cmcone{} face is defined on $M^2$, $G$ and $Q$ are as well,
 so $ds_\#^2$ is well-defined on $M^2$, and is called the {\em lift metric}.
 It is nothing but the {\em dual metric\/}
 of the \cmcone{} surface $\pi_H\circ F$ in $H^3$,
 see \cite{UY5}.
\end{remark}

\subsection*{Completeness}
We now define two different notions of
completeness for \cmcone{} faces as follows: 

\begin{definition}\label{def:c}
 We say a \cmcone{} face $f\colon M^2 \to S^3_1$ is
 {\it complete\/} if there exists a symmetric
 $2$-tensor field $T$ 
 which vanishes outside a compact subset $C\subset M^2$
 such that the sum $ T+ ds^2 $ is a complete Riemannian metric on
 $M^2$. 
\end{definition}

See \cite{F}, with similar definitions in
 \cite{KUY} for flat fronts in $H^3$ and in \cite{UY1} for maxfaces.
\

\begin{definition}\label{def:wc}
 We say that $f$ is {\it weakly complete\/} if
 it is congruent to an $S^3_1$-horosphere  or if
 the lift metric \eqref{eq:lift-metric}
 is a complete Riemannian metric on $M^2$.
\end{definition}

 Here, the 
 {\em $S^3_1$-horosphere} is the totally umbilic \cmcone{} surface,
 which is also the only complete \cmcone{} immersed surface
 (see Remark~\ref{rem:complete}).
 It has the Weierstrass data $g=c=\mbox{constant}$
 $(|c|\ne 1)$ and $\omega=dz$.
 The metric $ds^2_\#$ of an $S^3_1$-horosphere cannot be defined
 by \eqref{eq:lift-metric} as $G$ is constant and 
 $Q$ is identically $0$, 
 but can still be defined as the metric induced by $F^{-1}$, and is a
 complete flat metric on $\C$.

\begin{definition}\label{def:ft}
 We say that $f$ is of {\em finite type\/}
 if there exists a compact set $C$ of $M^2$ such that
 the first fundamental form $ds^2$ is positive definite
 and has finite total (absolute) curvature on $M^2\setminus C$.
\end{definition}

 Let $f\colon{}M^2\to S^3_1$ be a \cmcone{} face of
 finite topology, that is, $M^2$ is diffeomorphic to 
 a compact Riemann surface $\overline{M}^2$ with a 
 finite number of points $\{p_1,\dots,p_n\}\subset\overline{M}^2$ excluded. 
 We can take a punctured neighborhood $\Delta_j^*$ of $p_j$
 which is biholomorphic to either the punctured unit disk 
 $\Delta^*=\{z\in\C\,;\,0<|z|<1\}$ or
 an annular domain, and $p_j$ is called a {\it puncture-type end\/} or an
 {\it annular end},
 respectively.
\begin{proposition}\label{prop:complete}
 Let $f:M^2\to S^3_1$ be a \cmcone{} face.
 If $f$ is complete, then
 \begin{enumerate}
\item\label{item:complete:0}
 the singular set of $f$ is compact, 
  \item\label{item:complete:1}
        $f$ is weakly complete,
  \item\label{item:complete:2}
        $M^2$ has finite topology and each end
        is of puncture-type.
 \end{enumerate}
\end{proposition}

\begin{proof}
 \ref{item:complete:0} is obvious.  
 If $f$ is totally umbilic, it is congruent
 to an $S^3_1$-horosphere and the assertion is obvious.  
 So we assume the Hopf differential $Q$
 does not vanish identically.
 Since the Gaussian curvature of $f$ is nonnegative,
 completeness implies \ref{item:complete:2}
 by the appendix of \cite{UY1}.
 So we shall now prove that completeness implies weak completeness:
 Fix an end $p_j$ of $f$.  By an appropriate choice of a coordinate $z$,
 the restriction of $f$ to a neighborhood of
 $p_j$ is $f_j:\Delta^*\to S^3_1$.
 We denote by $d\hat s^2$ the induced metric of the corresponding 
 \cmcone{} surface 
 $\hat f_j=FF^*\colon{}\widetilde{\Delta}^*\to H^3$ into hyperbolic
 $3$-space.
 Take a path $\gamma\colon{}[0,1)\to \Delta^*$ such that 
 $\gamma(t)\to 0$ as $t\to 1$.
 Then by \eqref{eq:dshat^2} and \eqref{eq:forms}, 
 $d\hat s^2 \geq ds^2$ holds, and hence
 completeness of $f$
 implies that each lift
 $\tilde\gamma\colon{}[0,1)\to\widetilde{\Delta}^*$
 of $\gamma$ has infinite length with respect to $d\hat s^2$.
 Here, $d\hat s^2$ and 
 $ds^2_{\#}=(1+|G|^2)^2 |Q/dG|^2$ are the pull-backs of the 
 Hermitian metric of $\SL_2\C$ by $F$ and $F^{-1}$, respectively.
 Yu \cite{Yu} showed that completeness of these two metrics are
 equivalent.
 Hence, $\tilde\gamma$ has 
 infinite length with respect to the metric $ds^2_{\#}$.
 Since $ds^2_{\#}$ is well-defined on $\Delta^*$, 
 $\gamma$ also has infinite length with respect to  $ds^2_{\#}$, 
 that is, the metric $ds^2_{\#}$ on $\Delta^*$ is complete at $0$.
 Thus, $f_j$ is a weakly complete end.
\end{proof}

For further properties of complete ends, see
Theorems~\ref{thm:completeend} and \ref{Thm:labeledbyYang3}.
\begin{remark}
 Our definition of weak completeness  of \cmcone{} faces  is
 somewhat more technical than that of maxfaces \cite{UY1}, 
 but
 it is the correctly corresponding concept in $S_1^3$: for
 data $(g,\omega)$, weak completeness of the associated maxface 
 in $\R_1^3$ is equivalent to that of the 
 \cmcone{} face in $S^3_1$.
\end{remark}

\begin{remark}
The \cmcone{} trinoids in $S^3_1$ constructed in \cite{LY} 
are all weakly complete (sometimes complete as well) and
all ends are $g$-regular, see Section~\ref{sec:ends}.
\end{remark}

\begin{remark}
The Hopf differential $Q$ of a complete \cmcone{} face
$f:M^2\to S^3_1$ is meromorphic on its compactification
$\overline M^2$, even without assuming that all ends of $f$
are regular.  See Appendix~\ref{app:Q-mero}.
It should be remarked that for \cmcone{} surfaces in hyperbolic $3$-space,
finiteness of total curvature  
is needed to show the meromorphicity of $Q$ (see \cite{Br}).
\end{remark}

\subsection*{Monodromy of ends of \cmcone{} faces}
For any real number $t$, we set
\begin{equation}\label{eq:canonical-repr}
 \begin{aligned}
  \Lambda_e(t)&:=
  \begin{pmatrix}
     e^{it} & 0 \\
      0     & e^{-it}
  \end{pmatrix}, \\
  \Lambda_p(t) &:=
  \begin{pmatrix}
    1+i t & -it\\
    it & 1-it
  \end{pmatrix},  \\
  \Lambda_h(t) &:=
  \begin{pmatrix}
    \cosh t & \sinh t\\
    \sinh t & \cosh t
  \end{pmatrix}.
 \end{aligned}
\end{equation}
A matrix in $\SU_{1,1}$ is called
\begin{enumerate}
 \item {\em elliptic\/} if it is conjugate 
       to $\Lambda_e(t)$ ($t\in(-\pi,\pi]$)
       in $\SU_{1,1}$,
 \item {\em parabolic\/} if it is conjugate
       to $\pm\Lambda_p(t)$ ($t\in\R\setminus\{0\}$) in $\SU_{1,1}$, and
 \item {\em hyperbolic\/}
       if it is conjugate
       to $\pm\Lambda_h(t)$ ($t>0$) in $\SU_{1,1}$.
\end{enumerate}
Any matrix in $\SU_{1,1}$ is of one of
these three types, see Appendix~\ref{app:conj}.
Note that the parabolic matrices $\Lambda_p(t_1)$ and $\Lambda_p(t_2)$
are conjugate in $\SU_{1,1}$ if and only if $t_1t_2>0$.
Though the set of conjugate classes of parabolic matrices 
is fully represented by
$\{\pm\Lambda_p(\pm 1)\}$,
we may use various values of $t$ in this paper for the sake of
simplicity.  

Let $f\colon M^2 \to S^3_1$ be a weakly complete \cmcone{} face of
finite topology,
where $M^2$ is diffeomorphic to a compact Riemann surface
$\overline{M}^2$ with
finitely many punctures $\{p_1 , \dots , p_n\}$.
Any puncture $p_j$, or occasionally
a small neighborhood $U_j$ of $p_j$, is called an {\em end\/} of $f$.

An end is called {\it elliptic}, {\it parabolic} or {\it hyperbolic\/}
when the monodromy matrix
$\tilde\rho(\tau)\in\SU_{1,1}$ 
is elliptic, parabolic or hyperbolic, respectively,
where $\tilde{\rho}$ is as in Remark~\ref{rem:ambiguity}
and $\tau\in\pi_1(M^2)$ is the deck transformation 
corresponding to the counterclockwise loop about $p_j$.

\subsection*{The Schwarzian derivative}
Let $(U,z)$ be a local complex coordinate of a Riemann surface $M^2$,
and $h(z)$ a meromorphic function on $U$.  Then
\[
   S_z(h) :=
    \left(\frac{h''}{h'}\right)'-
    \frac{1}{2}
    \left(\frac{h''}{h'} \right)^2
    \qquad \left('=\frac{d}{dz}\right)
\]
is the {\em Schwarzian derivative\/} of
$h$ with respect to the coordinate $z$.

If $h(z)=a+b (z-p)^m+o\bigl((z-p)^m\bigr)$ at $z=p$ $(b\ne 0)$,
where $o\bigl((z-p)^m\bigr)$ denotes higher order terms,
then the positive integer 
$m$ is called the {\em {\rm (}ramification{\rm )} order} of $h(z)$,
and we have
\begin{equation}\label{eq:s_ord}
 S_z(h)=\frac{1}{(z-p)^2} \left( \frac{1-m^2}{2}+ o(1)\right).
\end{equation}
We write $S(h)=S_z(h)\,dz^2$, which we also call the 
{\em Schwarzian derivative.
The Schwarzian derivative depends on
the choice of local coordinates, but the difference
does not, that is, 
$S(h_1)-S(h_2)$ is a well-defined holomorphic $2$-differential.}

The Schwarzian derivative is invariant under M\"obius transformations:
$S(h)=S(A\star h)$ holds for $A\in \SL_2\C$, where  $\star$
denotes the M\"obius transformation as in \eqref{eq:second-moebius}.
Conversely, if $S(h)=S(g)$, there exists an $A\in\SL_2\C$ such that
$g=A\star h$.  

Let $f\colon{}M^2\to S^3_1$ be a \cmcone{} face with the hyperbolic
Gauss map $G$, a secondary Gauss map $g$ and the Hopf
differential $Q$.  Then 
\begin{equation}\label{eq:schwarz-Q}
  S(g)-S(G)= 2 Q.
\end{equation}

\begin{remark}\label{rem:complete}
 Here we give a  proof that 
 the only complete \cmcone{} immersion is the totally umbilic one,
 that is, the $S^3_1$-horosphere,
 which is simpler than the original proofs in \cite{Ak,R}. 
 (The proof is essentially the same as for the case of maximal surfaces
 in $\R^3_1$ given in \cite[Remark 1.2]{UY1}.)
 Let $f:M^2\to S^3_1$ be a complete \cmcone{} immersion.
 Without loss of generality, we may assume that $M^2$
 is both connected and simply connected. Then the Weierstrass
 data $(g,\omega)$ as in Theorem \ref{th:AA-rep} is 
 single-valued on $M^2$.
 Since $f$ has no singular points, we may assume that
 $|g|<1$ holds on $M^2$. 
 Since $(1-|g|^2)^2|\omega|^2<|\omega|^2$, the metric
 $|\omega|^2$ is a complete flat metric on $M^2$.
 Then the uniformization theorem yields that
 $M^2$ is bi-holomorphic to $\C$, and $g$ must be
 a constant function, which implies that
 the image of $f$ must be totally umbilic.
\end{remark}

\section{Monodromy of punctured hyperbolic metrics}
\label{sec:projective}
By Remark~\ref{rem:ambiguity}, the monodromy of a holomorphic null
immersion $F$ is elliptic, parabolic or hyperbolic if and only if the
monodromy of its secondary Gauss map $g$ is elliptic, parabolic or
hyperbolic, respectively. 
In this section, in an abstract setting, we give results needed for 
investigating the behavior of $g$ at a puncture-type end, in terms of the 
monodromy of $g$.  

\subsection*{Lifts of $\PSU_{1,1}$-projective 
connections on a punctured disk.}
Let
\[
   \Delta^*=\Delta\setminus\{0\},\qquad \mbox{where}\,\,
   \Delta:=\{z\in\C\,;\,|z|<1\},
\]
be the punctured unit disk and $P=p(z)dz^2$ 
a holomorphic $2$-differential on $\Delta^*$.
Then there exists a holomorphic developing map
$g_P\colon{}\widetilde\Delta^*\to\C\cup\{\infty\}$
such that
$S(g_P^{})=P$,
where $\widetilde\Delta^*$ is the universal cover of
$\Delta^*$.
For any other holomorphic function $h$ such that
$S(h)=P$, there exists an $A \in \SL_2\C$
so that $ A \star g_P^{}=h$.  Thus there exists a matrix
$T \in \PSL_2\C$ such that
\begin{equation}\label{eq:monodromy-punctured}
    g_P^{} \circ\tau ^{-1} = T\star g_P^{},
\end{equation}
where $\tau$ is the generator of $\pi_1(\Delta^*)$
corresponding to a counterclockwise loop about the origin.
We call $T$ the monodromy matrix of $g_P$. 
If there exists a $g_P^{}$ so that 
$T\in \PSU_{1,1}$, $P$ is called a
{\em $\PSU_{1,1}$-projective connection\/} on $\Delta^*$ and
$g_P^{}$ is called a $\PSU_{1,1}$-lift of $P$.
A $\PSU_{1,1}$-projective connection on $\Delta^*$ has a removable
singularity,
a pole or an essential singularity at $0$,
and is said to have a \emph{regular singularity}
at $0$ if it has at most a pole of order 2 at $0$.
(The general definition of projective connections is given in \cite{T}
and \cite{UYh}. 
There exist holomorphic $2$-differentials on $\Delta^*$ which are
not $\PSU_{1,1}$-projective connections.)
When $T \in \PSU_{1,1}$, 
it is conjugate to one of the matrices in
\eqref{eq:canonical-repr}.
The $\PSU_{1,1}$-projective connection $P$ is then called
{\it elliptic}, {\em parabolic} or {\em hyperbolic}
when $T$ is  {\it elliptic}, {\em parabolic} or {\em hyperbolic}, 
respectively.
This terminology is independent of the choice of $g_P^{}$.

By the property \eqref{eq:g-repr}, the Schwarzian derivative $S(g)$
of the secondary Gauss map $g$ of a \cmcone{} face is 
an example of a $\PSU_{1,1}$-projective connection.

Note that a $\PSU_{1,1}$-lift $g_P^{}$ has the $\PSU_{1,1}$ ambiguity
$g_P^{}\mapsto A\star g_P^{}$ for $A\in\PSU_{1,1}$.
The property that $|g_P^{}|>1$ (resp.~$|g_P^{}|<1$) is independent of
this ambiguity.
\begin{remark}\label{rem:inverse}
 Let $g_P$ be a $\PSU_{1,1}$-lift of a $\PSU_{1,1}$-projective
 connection $P$.
 Then 
 \[
    \frac{1}{g_P} = D\star g_P
    \qquad 
    \left(
     D:=\begin{pmatrix} 0 & i \\ i & 0 \end{pmatrix}\right)
 \]
 is also a $\PSU_{1,1}$-lift of $P$, because $DAD^{-1}\in \PSU_{1,1}$
 for any $A\in\PSU_{1,1}$.
 However, $D\not\in\PSU_{1,1}$, and 
 one can show that there is no matrix $B\in\PSU_{1,1}$ such that
 $1/g_P=B\star g_P$, that is, 
 $1/g_P$ is not $\PSU_{1,1}$-equivalent to $g_P$.
\end{remark}

In the rest of this article, as well as in the following
proposition,
we use
\begin{equation}\label{matrixR}
  R := \frac{1}{2}\begin{pmatrix}
          1 & \hphantom{-}1 \\ i & -i
         \end{pmatrix},
\end{equation}
which is motivated by an isomorphism between 
$\SL_2 \R$ and $\SU_{1,1}$.
See Appendix~\ref{app:conj}.
\begin{proposition}\label{N}
 Let $P$  
 be a  $\PSU_{1,1}$-projective connection on $\Delta^*$.
 Then the following assertions hold{\rm:}
 \begin{enumerate}
  \item\label{item:elliptic}
    Suppose that $P$ is elliptic. Then,
   \begin{enumerate}       
    \item\label{item:elliptic:1}
          there exist a real number
          $\mu$ 
          and a single-valued meromorphic function $h(z)$ on 
          $\Delta^*$ such that
          \[  
                g(z) :=z^{\mu} h(z)
          \] 
          is a $\PSU_{1,1}$-lift of $P$.
    \item\label{item:elliptic:2} 
          $P$ has a regular singularity
          at $z=0$ if and only if  $h(z)$ has at most a pole at $z=0$.
   \end{enumerate}
  \item\label{item:parabolic}
    Suppose that $P$ is parabolic and 
    take an arbitrary positive number $t$.
    Then,
    \begin{enumerate}
     \item\label{item:parabolic:1}
           for each $\varepsilon\in\{-1,1\}$,
           there exists a 
           single-valued meromorphic function $h(z)$ on 
           $\Delta^*$ such that
           \[
              g(z) := R^{-1}\star 
                      \left(
                         h(z)-\frac{\varepsilon t}{\pi i}\log z
                      \right)
           \]
           is a $\PSU_{1,1}$-lift of $P$.
     \item\label{item:parabolic:2}
           The function $h(z)$ has at most a pole at $z=0$ if and only
           if $P$ has a pole of order exactly $2$ at $z=0$.
     \item\label{item:parabolic:3}
           $h(z)$ is holomorphic at $z=0$
           if and only if $P- dz^2/(2z^2)$ has at most a pole of order
           $1$ at $z=0$.
     \item\label{item:parabolic:4}
            When $h(z)$ is holomorphic at $z=0$, $|g(z)|>1$ {\rm
           (}resp.~$|g(z)|<1${\rm )} 
           holds
           for sufficiently small $|z|$ 
           if and only if $\varepsilon=+1$ 
           {\rm (}resp.~$\varepsilon=-1${\rm )}.
    \end{enumerate}
  \item\label{item:hyperbolic}
    Suppose that $P$ is hyperbolic.
    Then,
    \begin{enumerate}
     \item\label{item:hyperbolic:1}
           there exist a positive number $\mu$
           and a
           single-valued meromorphic function $h(z)$ on 
           $\Delta^*$ such that
           \[
                 g(z) := R^{-1}\star\left(z^{i\mu} h(z)\right)
           \]
              is a $\PSU_{1,1}$-lift of $P$.
     \item\label{item:hyperbolic:2}
           $h(z)$ has at most a pole at $z=0$ if and only
           if $P$ has a pole of order exactly $2$ at $z=0$.
    \end{enumerate}
 \end{enumerate}
\end{proposition}
\begin{remark}
 In the statements of Proposition~\ref{N}, 
 the function  $z^{\mu}$ $(\mu\in\C)$ is defined by
 \[
    z^{\mu} := \exp(\mu\log z),
 \]
 where  $\log z$ is considered as 
 a function defined
 on the universal cover $\widetilde\Delta^*$ of $\Delta^*$.
\end{remark}

To prove this, we consider
the following ordinary differential equation
\begin{equation}\label{ode}
    X''+\frac{1}{2} p(z) X=0 \qquad 
     \left(~'=\frac{d}{dz},~P=p(z)\,dz^2\right).
\end{equation}
If we assume $P(z)$ has a regular singularity at $z=0$, then
$p(z) = \alpha z^{-2}\bigl(1+o(1)\bigr)$ for some 
$\alpha \in \C$
and \eqref{ode} has the fundamental system of solutions 
\begin{equation}\label{eq:fundamental}
 \begin{aligned}
  X_1(z)&=z^{\mu_1}\xi_1(z),\\
  X_2(z)&=z^{\mu_2}\xi_2(z) + k \log z X_1
 \end{aligned}
 \qquad(\Re \mu_1\geq\Re\mu_2 ),
\end{equation}
where $\xi_j(z)$ $(j=1,2)$ are
holomorphic functions on $\Delta=\{|z|<1\}$ such that
$\xi_j(0)\ne 0$ $(j=1,2)$.  The constant $k\in \C$ is
called the {\it log-term coefficient\/} and
$\mu_1,\mu_2$ are the solutions of the indicial equation
\begin{equation}\label{eq:indicial}
   t(t-1)+\frac{\alpha}2=0.
\end{equation}
If $\mu_1-\mu_2 \not\in \Z$, then $k$ vanishes.
(See \cite{CL} or the appendix of \cite{RUY-hiroshima}).
The following lemma is easy to show:
\begin{lemma}\label{lem:ode}
 In the above setting,
 $S(g_0)=P$ if $g_0:={X_2}/{X_1}$.
\end{lemma}
\begin{proof}[Proof of Proposition~\ref{N}]
 Take the matrix $T$ as in \eqref{eq:monodromy-punctured}.

 We first prove the elliptic case.
 Since $P$ is elliptic, there exist a
 $t \in \R$ and an $A\in\SU_{1,1}$ such that
 $ATA^{-1} = \Lambda_e(t)$.
 So $(A\star g)\circ\tau^{-1}=e^{2it}(A\star g)$,
 and $h(z):=z^{t/\pi}\bigl(A\star g(z)\bigr)$ is single-valued
 on $\Delta^*$, proving the first part of \ref{item:elliptic}.
 If the origin $0$ is at most a pole of $h$,
 a direct calculation shows that $P$ has a regular singularity.
 To show the converse,
 we set
 $g_0:=X_2/X_1$,
 with $\{X_1,X_2\}$ as in \eqref{eq:fundamental}.
 Then by Lemma~\ref{lem:ode} we have $S(g_0)=P$.
 The monodromy matrix $\pm T_0$ of $g_0$ is conjugate
 to 
\begin{equation}\label{eq:monodromy}
\begin{cases}
 \begin{pmatrix}
  e^{\pi i (\mu_1-\mu_2)}  &  0 \\ 0 & e^{\pi i (\mu_2-\mu_1)} 
 \end{pmatrix}
 & \qquad\mbox{(if $k=0$)}, \\[12pt]
 \begin{pmatrix}
  1  &  -2\pi i k \\ 0 & 1
 \end{pmatrix}
 & \qquad\mbox{(if $k\ne 0$)}.
\end{cases}
\end{equation}
 Since $P$ is elliptic, the log-term coefficient $k=0$
 and $\mu_2-\mu_1 \in \R$.  Thus
 \[
    g_0(z)=z^{\mu} \frac{\xi_2(z)}{\xi_1(z)} \qquad
        (\mu:=\mu_2-\mu_1).
 \]
 Since $S(A\star g)=S(g_0)$, there exists a $B\in \SL_2\C$ so that
 $A\star g=B\star g_0$.  
 Then
 \[
  \Lambda_e(t) \star (A\star g)
  =(A\star g)\circ \tau^{-1}
  =B\star (g_0\circ \tau^{-1})
   =B
    \Lambda_e(-\pi \mu)
     B^{-1}\star (A\star g),
 \]
 so $\Lambda_e(t) = \pm B \Lambda_e(-\pi \mu) B^{-1}$.
 If $t \equiv 0 \pmod{\pi}$, then $A\star g$ is
 meromorphic, proving \ref{item:elliptic}.
 Otherwise, 
 \[
     B = \begin{pmatrix}  
          c & 0 \\ 0 & c^{-1}\end{pmatrix}
     \qquad \text{or}\qquad
     \begin{pmatrix}
            0  & c \\
          -c^{-1} & 0 
     \end{pmatrix}
 \]
 for some $c\in \C\setminus\{0\}$,
 and \ref{item:elliptic} follows from
 \[
    A\star g(z) = c^2 z^{\mu} \frac{\xi_2(z)}{\xi_1(z)},\qquad
     \text{or}\qquad
    -c^{2} z^{-\mu} \frac{\xi_1(z)}{\xi_2(z)},
 \]
 respectively.

 Next, we assume $P$ is parabolic and take a positive number $t$
 and $\varepsilon\in\{-1,1\}$.
 Then by Theorem~\ref{thm:canonical-su} and 
 Remark~\ref{rem:parab-conj} in Appendix~\ref{app:conj},
 there exists a matrix $A\in\SU_{1,1}$ such that
 $ATA^{-1}$ is one of 
 $\Lambda_p(\varepsilon t)$, 
 $-\Lambda_p(\varepsilon t)$, $\Lambda_p(-\varepsilon t)$, 
 $-\Lambda_p(-\varepsilon t)$.
 Note that $\Lambda_p(\varepsilon t)$ and $\Lambda_p(-\varepsilon t)$
 are not conjugate in $\PSU_{1,1}$.
 Replacing $g$ with $1/g$ 
 if $ATA^{-1}=\pm\Lambda_{p}(-\varepsilon t)$
 (see Remark~\ref{rem:inverse}), 
 we can choose a $\PSU_{1,1}$-lift $g$ such that
 \[
   ATA^{-1}   = \pm\Lambda_p(\varepsilon t).
 \]
 Then, 
 $(A\star g)\circ \tau^{-1}= \Lambda_p(\varepsilon t)\star 
 (A\star g)$. 
 Here the $\pm$-ambiguity of $ATA^{-1}$ does not affect 
 the $\star$-action.
 Thus, 
 \[
    \bigl((RA)\star g\bigr)\circ \tau^{-1}=(RA)\star g +2\varepsilon t, 
 \qquad\text{since}\quad
    R \Lambda_p(\varepsilon t) =
    \begin{pmatrix} 1 & 2\varepsilon t \\ 0 & 1 \end{pmatrix}R.
 \]
 Hence
 $h(z):=(RA)\star g+\bigl(\varepsilon t/(\pi i)\bigr) \log z $
 is  a single-valued meromorphic function on $\Delta^*$,
 proving the first part of \ref{item:parabolic}.
 If $h(z)$ has at most a pole at $z=0$, then a direct computation
 shows that $P$ has a pole of order exactly $2$. 
 Therefore, it suffices to show that $h(z)$ has at most a pole
 at $z=0$ when $P$ has a regular singularity. We now show this:

 We set $g_0:=X_2/X_1$.
 Since $P$ is parabolic, \eqref{eq:monodromy} yields that
 the log-term coefficient $k \neq 0$
 and $\mu:=\mu_2-\mu_1 \in \Z$ is non-positive. 
 Hence
 \[
    g_0=z^{\mu} \frac{\xi_2(z)}{\xi_1(z)}+k \log z.
 \]
 Here $z^{\mu} \xi_2(z) / \xi_1(z)$
 is single-valued on $\Delta^*$ and has at most a pole at $z=0$.
 Take a matrix $B\in \SL_2\C$ such that $(RA)\star g=B\star g_0$.  Then
 \begin{align*}
  \begin{pmatrix}
    1 & 2\varepsilon t \\ 0 & 1
  \end{pmatrix}&\star (RA\star g)
  =(RA\star g)\circ \tau^{-1}
     =B\star (g_0\circ \tau^{-1})\\
  &=B\begin{pmatrix}
      1 & -2 \pi i k \\ 0 & 1
     \end{pmatrix}\star g_0
  =B
  \begin{pmatrix}
   1 & -2 \pi i k \\ 0 & 1
  \end{pmatrix} B^{-1}
  \star (RA\star g).
 \end{align*}
 Replacing $X_2$ by $(i \varepsilon t/(\pi k)) X_2$ and renaming
 $(i\varepsilon t/(\pi k)) \xi_2$ to $\xi_2$, $-2\pi i k$ becomes 
 $2 \varepsilon t$, 
and we have
 \[
   \begin{pmatrix}
    1 & 2\varepsilon t \\ 0 & 1
    \end{pmatrix}
    =B  \begin{pmatrix}1 & 2\varepsilon t \\ 0 & 1\end{pmatrix} B^{-1}, 
 \]
 and there is no $\pm$-ambiguity in the above equation, as
 the eigenvalues of the left hand matrix
 must have the same sign as those of the right hand matrix.
 Thus we can choose
 \[
     B=
      \begin{pmatrix}
       1 & c \\ 0 & 1
      \end{pmatrix}
 \]
 for some $c\in \C$,
 which proves the second part of \ref{item:parabolic}.
 It is easy to see that
 $h(z)$ is holomorphic at $z=0$ if and only if $\mu_2=\mu_1$, 
 that is, $\alpha$ in \eqref{eq:indicial} is $1/2$,
 which proves the third part of \ref{item:parabolic}.
 Assume $h$ is holomorphic on $\Delta$.
 Since the M\"obius transformation $z\mapsto R \star z$
 maps the disk $\Delta$ onto the upper-half plane $\{z;\Im z>0\}$,
 the condition $|g|>1$ (equivalently $|A\star g|>1$)
 is equivalent to $\Im (RA\star g)<0$ 
 for all $A\in\SU_{1,1}$.
 And since $|h|$ is bounded, this is equivalent to $\epsilon>0$.
 Thus we obtain the last part of \ref{item:parabolic}.

 Next, we assume $P$ is hyperbolic.
 By Theorem~\ref{thm:canonical-su} in Appendix~\ref{app:conj},
 there are a matrix $A\in\SU_{1,1}$ and $t >0$
 such that
 $ATA^{-1}=\Lambda_h(t)$ or $-\Lambda_h(t)$.
 Then $(A\star g)\circ \tau^{-1}=\Lambda_h(t)\star (A\star g)$,
 which implies $\bigl((RA)\star g\bigr) \circ \tau^{-1}=e^{2t}(RA)\star g$.
 So $h(z):=z^{-i t/\pi}(RA)\star g$
 is a single-valued meromorphic function in $\Delta^*$.
 This proves the first part of \ref{item:hyperbolic}.
 To prove the second part of \ref{item:hyperbolic},
 analogous to the parabolic case, we only need to prove 
 one direction. Suppose that $P$ has a regular singularity.
 We set $g_0:=X_2/X_1$.  Since $P$ is hyperbolic, 
 \eqref{eq:monodromy} yields that
 $k=0$  and $\mu_2-\mu_1 \in i \R$.  Thus
 \[
   g_0(z)=z^{i\mu} \frac{\xi_2(z)}{\xi_1(z)} \; , \qquad
    \mu:=i (\mu_1-\mu_2).
 \]
 Exchanging $X_1$ and $X_2$ if necessary, we may assume $\mu>0$
 without loss of generality.
 Take a $B\in \SL_2\C$ such that $RA\star g=B\star g_0$.
 Then we have
 \[
  \begin{pmatrix} e^t & 0 \\ 0 & e^{-t}\end{pmatrix}
  \star (RA\star g)
  =B\star (g_0\circ \tau^{-1})
  =B 
  \begin{pmatrix} e^{\pi\mu} & 0 \\ 0 & e^{-\pi\mu}\end{pmatrix}
  B^{-1}\star (RA\star g),
 \]
 so
 \[
   \begin{pmatrix}
      e^t & 0 \\
      0 &e^{-t}
   \end{pmatrix}
    =B
    \begin{pmatrix}
       e^{\pi\mu} & 0 \\
       0 & e^{-\pi\mu}
    \end{pmatrix}
    B^{-1},
 \]
 that is, $t=\pm \pi\mu$. 
 As we have assumed that $t>0$ and $\mu>0$, 
 we have $t=\pi\mu$, and then
 $B$ must be diagonal.
 Hence
 \[
    RA\star g(z)=c^2 z^{it/\pi}
     \frac{\xi_2(z)}{\xi_1(z)}, \qquad
     B=\begin{pmatrix} c & 0 \\ 0 & c^{-1}\end{pmatrix}
 \]
 proving the assertion.  
\end{proof}

\subsection*{Monodromy of punctured hyperbolic metrics}

We consider a conformal metric
$d\sigma^2$ on $\Delta^*$ of constant Gaussian curvature $-1$,
called a {\it punctured hyperbolic metric}.
Then there exists a meromorphic function
$g:\widetilde{\Delta}^*\to \C\cup\{\infty\}\setminus\{|z|=1\}$
such that
\begin{equation}\label{eq:hyperbolic-metric}
   d\sigma^2 = \frac{4\,|dg|^2}{(1-|g|^2)^2},
\end{equation}
which is called the {\it developing map\/} of $d\sigma^2$.
Since $d\sigma^2$ is a well-defined hyperbolic metric on $\Delta^*$,
either $|g|<1$ or $|g|>1$ holds on $\Delta^*$.

We remark that,
for a \cmcone{} immersion $f\colon{}\Delta^*\to S^3_1$,
the metric $d\sigma^2$ as in \eqref{eq:hyp-metric} is 
an example of a hyperbolic metric,
and the secondary Gauss map is a developing map of it.

The developing map $g$ is not unique, 
and the set of all developing maps of
$d\sigma^2$ coincides  with 
\[
   \left\{A\star g; A\in \SU_{1,1}\right\}\cup
   \left\{A\star \frac{1}{g}=
          A\begin{pmatrix} 0 & i \\ i & 0 \end{pmatrix}\star g; A\in
            \SU_{1,1}\right\} \; . \]

Set
\begin{equation}\label{eq:metric-schwartz}
    S(d\sigma^2):=S(g)= \left(w_{zz}-\frac{(w_z)^2}{2}\right)dz^2,
\end{equation}
    where $d\sigma^2=e^w\, |dz|^2$, that is,
$w:=
 \log\left(4|g_z|^2/(1-|g|^2)^2\right)$.
We call the projective connection
$S(d\sigma^2)$ the {\it Schwarzian derivative\/} of $d\sigma^2$.
Since the metric $d\sigma^2$ is well-defined on $\Delta^*$,
the developing map $g$ is a $\PSU_{1,1}$-lift of
the  $\PSU_{1,1}$-projective  connection $S(d\sigma^2)$.

If 
$g$ is a developing map of $d\sigma^2= K_{ds^2} ds^2$ in
\eqref{eq:hyp-metric}, then
\[
   g\circ\tau^{-1} = T\star g\qquad\text{for some } T\in\PSU_{1,1}.
\]
If the matrix $T$ is {\em elliptic} (resp.~{\em parabolic}, 
{\em hyperbolic}),
the metric $d\sigma^2$ is said to have
elliptic (resp.~parabolic, hyperbolic) monodromy.

\begin{definition}\label{def:reg-sing}
 We say that a hyperbolic punctured metric
 $d\sigma^2$ has a {\em regular singularity\/}
 at the origin if $S(d\sigma^2)$ has
 a regular singularity at the origin,
 that is, it has at most a pole of order $2$.
\end{definition}

\begin{theorem}\label{thm:key}
 Any conformal hyperbolic metric on $\Delta^*$ has a 
 regular singularity at $z=0$.
\end{theorem}
\begin{proof}
 Let $g$ be a developing map of a conformal hyperbolic metric
 $d\sigma^2$ on $\Delta^*$.

 Suppose $d\sigma^2$ has elliptic monodromy.
 Since $d\sigma^2$ has no singular points on $\Delta^*$,
 $|g|<1$ or $|g|>1$ holds on $\Delta^*$.
 Since $1/g$ is also a developing map of  $d\sigma^2$, we may
 assume that $|g|<1$.
 By Proposition~\ref{N}, there exists a real number $\mu$ such that
 $h(z) := z^{-\mu} g(z)$ is a single-valued function on $\Delta^*$.
 Multiplying $h(z)$ by $z^k$ ($k\in\Z$),
 we may assume that  $-1 < \mu \le 0$ without loss of generality.
 Thus 
 \[
       |h(z)|=|z|^{-\mu} |g(z)|<|z|^{-\mu}<1,
 \]
 and $h(z)$ has more than two exceptional values, so has at most a pole at
 $z=0$, by the Great Picard theorem.
 Then by \ref{item:elliptic:2} in Proposition~\ref{N}, 
 $S(g)$ has a regular singularity at the origin.

 Suppose $d\sigma^2$ has parabolic monodromy.
 Applying Proposition~\ref{N} for the $\PSU_{1,1}$-projective 
 connection $S(d\sigma^2)$ with $\varepsilon=-1$ and $t=\pi$,
 we can take a $\PSU_{1,1}$-lift $g$ such that 
 \[
     h = \hat g + i \log z\qquad
   \left(
        \hat g(z) := R\star g(z)=\frac1{i}\frac{g(z)+1}{g(z)-1}
   \right)
 \]
 is a single-valued meromorphic function on $\Delta^*$,
 where $R$ is the matrix in \eqref{matrixR}.
 Since $d\sigma^2$ has no singular points on $\Delta^*$,
 $|g|>1$ or $|g|<1$ holds.
 In particular, because $z\mapsto R\star z$ maps
 the unit disk into the upper-half plane, 
 we have $\Im \hat g>0$ (resp.\ $\Im \hat g<0$)
 if $|g|<1$ (resp.\ $|g|>1$).
 Here, it holds that
 $$
|z \exp(ih)|=|\exp(i\hat g)|
                = \exp(-\Im\hat g).
$$
 Thus, 
 \[
 \begin{alignedat}{2}
    |z\exp(ih)|&=\exp(-\Im\hat g)<1\qquad &&(\text{if $|g|<1$}),\\
    \left|\frac{1}{z}\exp(-ih)\right|&
         =\exp(\Im\hat g)<1\qquad &&(\text{if $|g|>1$}).
 \end{alignedat}
 \]
 Thus by the Great Picard theorem, there exist an integer $m$ and
 a holomorphic function $\phi(z)$
 with $\phi(0)\ne 0$ such that $\exp(\pm ih(z))=z^m \phi(z)$,
 that is,
$$
\pm ih(z)=m\log z + \log \phi(z).
$$
 Since $h(z)$ is single-valued,
 $m$ must be  $0$. 
 Therefore, $h(z)$ can be extended to be holomorphic at $z=0$,
 and  then by \ref{item:parabolic:2} of Proposition~\ref{N},
 the origin is a regular singularity of $S(d\sigma^2)$.
\renewcommand{\qed}{\relax}
\end{proof}
To prove the hyperbolic case, we need the following 
\begin{fact}[Montel's theorem]
 \label{Montel}
 If a family of holomorphic functions $\{f_n\}_{n=1,2,3,\dots}$
 defined on a domain $D(\subset \C)$
 have two exceptional values in common, then 
 they are a normal family, 
 that is, there is a subsequence $\{f_{n_j}\}_{j=1,2,3,\dots}$ 
 such that either $\{f_{n_j}\}_{j=1,2,3,\dots}$ 
or $\{1/f_{n_j}\}_{j=1,2,3,\dots}$
 converges uniformly on every compact set in $D$.
\end{fact}
\begin{proof}[Proof of Theorem~\ref{thm:key}, continued]
 The proof for the hyperbolic case is parallel to the proofs of
 Propositions 4 and 5 in  \cite{Br}.

 Suppose $d\sigma^2$ has hyperbolic mono\-dromy.
 Again, we may assume that $|g|<1$ without loss of generality.
 By Proposition~\ref{N}, 
 again replacing $g$ by $A\star g$ for some $A\in \SU_{1,1}$
 if necessary, there exists a positive real number 
 $\mu$ such that
 $h(z) := z^{-i\mu}\bigl( R \star g(z)\bigr)$ 
 is a single-valued meromorphic
 function on $\Delta^*$. 
 The function
 \[
   \hat g(z):= R\star g(z)=
                          \frac{1}{i}\frac{g(z)+1}{g(z)-1}
 \]
 has neither zeros nor poles in $\Delta^*$, and $\Im \hat g>0$.
 We now define a set
 \[
   \Omega :=\left\{z\in \C\,;\, 0<|z|<1, 
                     |\arg z|<\frac{2\pi}{3}\right\}
 \]
 and analytic functions $\zeta$ and $f_n$ for $n=1,2,3,\dots$
 from 
 $\Omega$ to $\C$ by 
 $\zeta(z) := \hat{g}(z^2)$ and $f_n(z):=\hat g( z^2/ 2^{2n})$.
 Then $\{f_n\}_{n=1}^\infty$ is a family of holomorphic functions on $\Omega$.
 Since $\Im \hat g > 0$,
 we have $\Im f_n > 0$.
 Thus $\{f_n\}$ is a normal family by Montel's theorem. 
 Two possible cases arise.
 \paragraph{\bf Case 1:}  
  First we consider the case that
  a subsequence   $\{f_n\}$ converges to a holomorphic function
  uniformly on any compact subset of $\Omega$.
  Since 
  \begin{multline*}
    \Omega_l:= \left\{z\in \C\,;\, |z|=\frac{1}{2^{l}}, 
                                     |\arg z|\le 
                                     \frac{3}{5}\pi\right\}
               \cup
               \left\{z\in \C\,;\, |z|=\frac{1}{2^{l+1}}, 
                                     |\arg z|\le 
                                     \frac{3}{5}\pi\right\}
   \\
      \cup
       \left\{z\in \C\,;\, |\arg z|=\frac{3}{5}\pi, 
               \frac{1}{2^{l+1}} \leq |z| \leq \frac{1}{2^{l}} \right\}
  \end{multline*}
  for a positive integer $l \in \Z_+$ is a compact subset of
  $\Omega$,
  there exist a positive number 
  $M \in \R_+$ and an $n_0 \in \Z_+$ such that
  $|f_n(z)|<M$
  holds on $\Omega_1$ for $n\ge n_0$. 
  This implies that
   $|\zeta(z)|<M$ on $\Omega_{n+1}$ for $n\ge n_0$.
  Then by the maximum principle, we have
  \[
    |\zeta(z)|<M \quad\text{on}\quad
     \left\{z\in \C\,;\, \frac{1}{2^{n+1}} \leq |z| \leq \frac{1}{2^{n}},
            |\arg z|\le \frac{3}{5}\pi\right\}
  \]
  for each $n>n_0+1$. 
  Thus we have
  \[
    |\hat g(z^2)|=|\zeta(z)|<M\quad\text{on}\quad
    \left\{z\in \C\,;\, 0<|z| \leq \frac{1}{2^{n_0+1}}, 
    |\arg z|\le  \frac{3}{5}\pi\right\}.
  \]
  On the other hand, since 
  $e^{- \pi |\mu|} < |z^{-i\mu}| < e^{ \pi |\mu|}$ for $|z|<1$ 
  and $|\arg z|<\pi$,   the function $h(z)$ is bounded in a 
  punctured neighborhood of $z=0$ and
  has a removable singularity there.
 \paragraph{\bf Case 2:}
  Next we consider the case that a subsequence 
  $\{1/f_n\}$ 
  converges to a holomorphic function $f$.
  Then we can conclude that $1/h(z)$
  is bounded on $\Delta^*$.
  In this case $h(z)$ has at most a pole at the origin.

  In both cases,
  $ S(g) = S\left( R^{-1} \star ( z^{i\mu} h(z)) \right)$ 
  has at most a pole of order two at $z=0$.
\end{proof}
\begin{remark}\label{rem:d-sigma}
 In Corollary~\ref{h-monodromy}, 
 we shall show that in fact the monodromy of $d\sigma^2$ can never be
 hyperbolic.
\end{remark}

\section{Intrinsic behavior of regular ends}
\label{sec:ends}
Let $f\colon{}M^2\to S^3_1$ be a weakly complete \cmcone{} face of
finite topology, and let $\overline{M}^2$ be a compact Riemann surface
such that $M^2$ is diffeomorphic to
$\overline{M}^2\setminus\{p_1,\dots,p_n\}$.

\begin{definition}\label{def:reg}
A puncture-type end $p_j$ of $f$ is called {\it regular\/} if the
hyperbolic Gauss map $G$ has at most a pole at $p_j$. 
\end{definition}

\begin{definition}\label{def:g-reg}
On the other hand, we say
a puncture-type end $p_j$ is  {\it $g$-regular}
if the
Schwarzian derivative $S(g)$ of the secondary Gauss map $g$ has at
most a pole of order $2$ at $p_j$,
that is, the pseudometric $d\sigma^2:=4|dg|^2/(1-|g|^2)^2$
has a regular singularity at $p_j$
(cf. Definition \ref{def:reg-sing}).
\end{definition}

When $g$ is single-valued, $g$-regularity implies that $g$ has at most
a pole at the end.
When the Hopf-differential has at most a pole of order $2$,
regularity and $g$-regularity are equivalent, 
by \eqref{eq:schwarz-Q}.

Theorem~\ref{thm:key}
can now be stated in terms of \cmcone{} faces as follows:
\begin{lemma}\label{lem:key2}
  All ends of a complete \cmcone{}  face are $g$-regular.
\end{lemma}
\begin{proof}
 By Proposition~\ref{prop:complete}, all ends are of puncture-type.
 So we can set $M^2=\overline{M}^2\setminus\{p_1,\dots,p_n\}$,
 where $\overline{M}^2$ is a compact Riemann surface.
 Let $(g,\omega)$ be a Weierstrass data for $f$.
 Since the singular set is compact, the metric $d\sigma^2$
 as in \eqref{eq:hyp-metric}
 is a punctured hyperbolic metric in a
 punctured neighborhood of $p_j$.
 Then $d\sigma^2$ has a regular singularity by Theorem~\ref{thm:key}, 
 and hence  $f$ is $g$-regular at $p_j$.
\end{proof}

\begin{figure}
\begin{center}
\includegraphics[width=.30\linewidth]{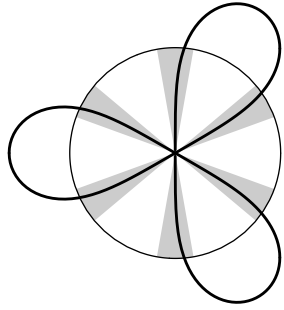}
\end{center}
\caption{The thick curve is the singular set of Example 
         \ref{ex:integral-elliptic-ends} for the case $m=3$. 
         The shaded parts indicate the set $S(3,\varepsilon,\pi/6)$.}
\label{fg:sing-e}
\end{figure}

\begin{definition}\label{def:labeledbyYang2}
 An elliptic end of a \cmcone{} face is
 {\it integral\/} if the monodromy of the secondary
 Gauss map is the identity, and
 {\it non-integral} otherwise.
\end{definition}
\begin{lemmaE1}
 Let $f\colon{}\Delta^*\to S^3_1$ be 
 a $g$-regular non-integral elliptic end.
 Then the singular set does not accumulate at the end $0$.
\end{lemmaE1}
\begin{proof}
 One can take the secondary Gauss map $g$ to be 
 $g(z) = z^{\mu} h(z)$ on a neighborhood of the end, 
 where $\mu\in \R\setminus \Z$ and $h(z)$ is holomorphic at the end $z=0$
 with $h(0)\ne 0$.
 Since $\mu\neq 0$, 
 the singular set $\{|g|=1\}$ cannot accumulate at the origin.
\end{proof}
On the other hand, an integral elliptic end 
might or might not be complete:
\begin{example}\label{ex:integral-elliptic-ends}
 For non-zero integers $m$ and $n$
 with $|m|\neq |n|$,
 we set $g=1-z^m$ and $G=z^n$.  
 Setting $Q=(S(g)-S(G))/2$ and $\omega =Q/dg$, we see that \eqref{eq:dshat^2} 
 gives a Riemannian metric on $\C\setminus\{0\}$.  So using Small's formula 
 \eqref{eq:small}, we have a \cmcone{} face with integral elliptic ends at 
 $z=0, \infty$.  
 The singular set is 
 $$
\{z\in\C\setminus\{0\}\,;\,|z|^{2m}-2\Re(z^m)=0\}
$$ 
 (see Figure \ref{fg:sing-e} for the case $m=3$).  
 Thus the singular set accumulates at $z=0$ but not at $z=\infty$.  
 Thus $z=\infty$ is a complete integral elliptic end, 
 but $z=0$ is an incomplete integral elliptic end. 
\end{example}

\begin{figure}
\begin{center}
\begin{tabular}{c@{\hspace{.1\linewidth}}c} 
\includegraphics[width=.30\linewidth]{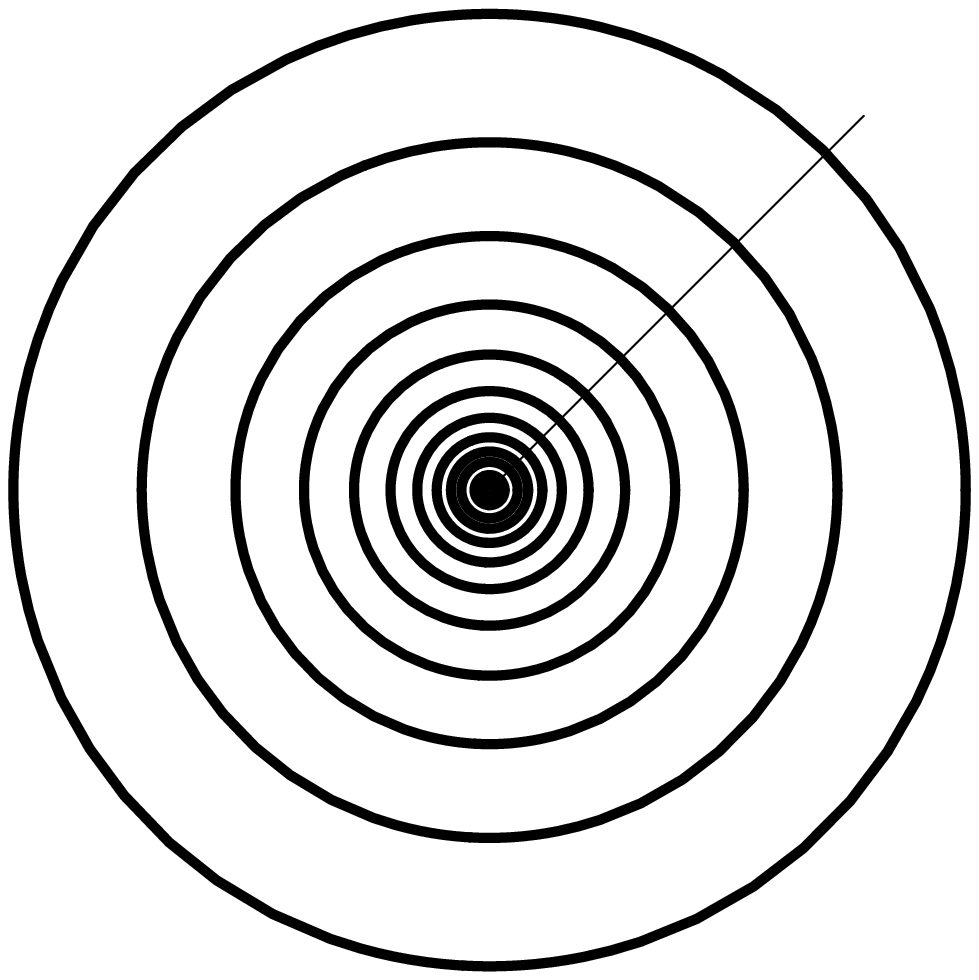} &
\includegraphics[width=.30\linewidth]{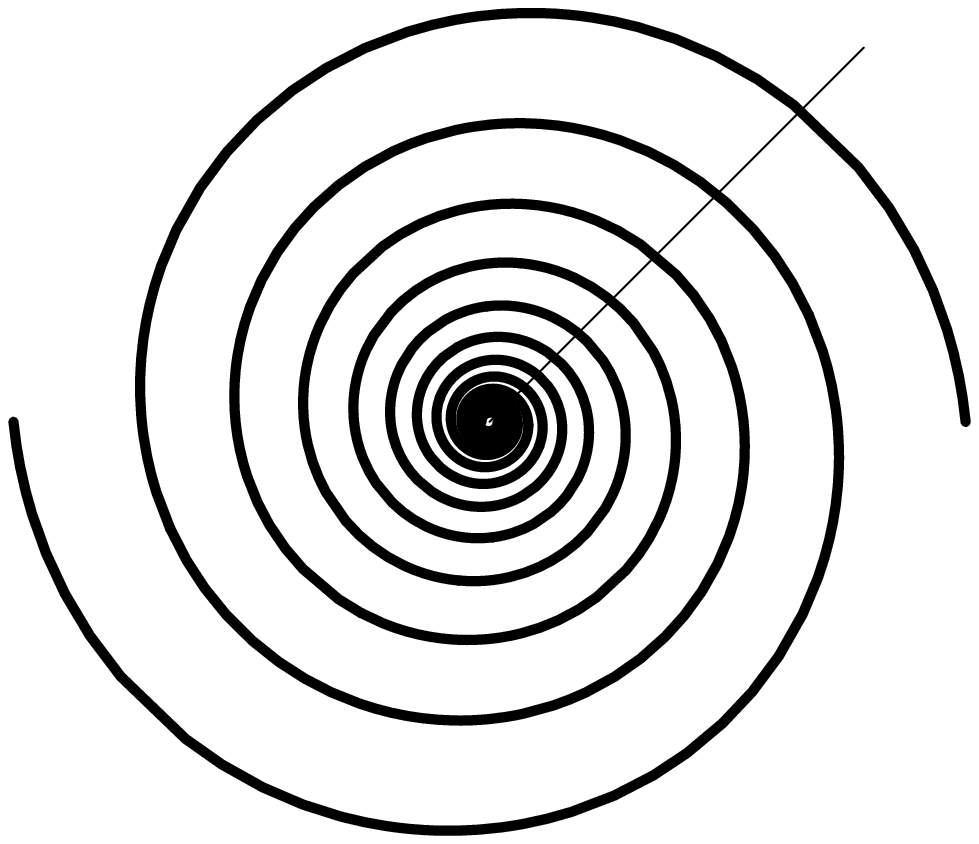} \\
$\hat g=w^{10i}$ & $\hat g=w^{1+10i}$
\end{tabular}
\end{center}
\caption{The thick curves indicate the singular sets of 
         hyperbolic ends with secondary Gauss map $g=R^{-1}\star\hat g$, 
         as discussed in Lemma H. The thin lines indicate rays 
         in $\Delta^*$ emanating from the origin.
         Here we see that the singular sets intersect rays emanating 
         from the origin infinitely many times.}
\label{fg:sing-h}
\end{figure}

 To state the behavior of an incomplete (integral) elliptic end, we
 introduce a notation:
 For a positive integer $m$, an $\varepsilon \in (0,\pi/(2m))$ 
 and a $\delta\in[0,\pi/m]$,
 we define the open subset 
 (which is a union of sectors, see Figure \ref{fg:sing-e})
 \[
    S(m,\varepsilon,\delta):=\bigcup_{k=0}^{2m-1} \left\{
        z\in \C\setminus\{0\}  \,;\, 
                \frac{k}{m}\pi+\delta-\varepsilon <\arg z<
                      \frac{k}{m}\pi+\delta+\varepsilon \right\}  .
 \]
\begin{lemmaE2}
 Suppose $f:\Delta^*\to S^3_1$ is a $g$-regular integral elliptic end.
 If the singular set accumulates at the end,
 then there are an $m \in \Z_+$ and a $\delta \in \R$ such that,
 for any $\varepsilon>0$,
 there exists an $r>0$ so that
 the singular set of $f$ in $\{z; 0<|z|<r\}$ lies in
 $S(m,\varepsilon,\delta)$.
\end{lemmaE2}
\begin{proof}
 The assertion of the lemma  does not depend on a choice of the complex
 coordinate at the origin.

 Since the singular set accumulates at $0$,
 we have $|g(0)|=1$.
 Then by Proposition~\ref{N} \ref{item:elliptic}, 
 $g(z)$ is holomorphic at $z=0$.
 Moreover, we may set $g(0)=1$.
 Then $\varphi(z):=\log g(z)$ is well-defined 
 on a neighborhood of $z=0$ and $\varphi(0)=0$.
 Using the Weierstrass preparation theorem,
 we may further assume without loss of generality that
 $\varphi(z)=z^m$ for some positive integer $m$.
 Here, $|g(z)|=1$ is equivalent to $\Re \varphi(z)=0$.
 Thus, the singular set is expressed as $\{\cos m\theta=0\}$,
 where $z=re^{i\theta}$.
\end{proof}
\begin{definition}\label{def:first-second}
 A parabolic end of a \cmcone{} face is of 
 the {\it first kind\/} if
 \[
          S(d\sigma^2)-\frac{dz^2}{2z^2}=
          S(g)-\frac{dz^2}{2z^2}=S(G)+2Q-\frac{dz^2}{2z^2}
 \]
 has at most a pole of order $1$. 
 Otherwise,
 it is of the {\it second kind}.
\end{definition}
\begin{lemmaP}
 Let $f:\Delta^*\to S^3_1$ be a $g$-regular parabolic end.
 If the end is of the first kind,
 the singular set does not accumulate at the end.
 If the end is of the second kind,
 then the singular set does accumulate at the end.
 In this case, there exist an $m \in \Z_+$
 and a $\delta$ $(\delta\in[0,\pi/m])$ such that, 
 for all $\varepsilon > 0$,
 there exists an $r>0$ so that
 the singular set of $f$ in $\{ z; 0<|z|<r\}$ lies in
 $S(m,\varepsilon,\delta)$.
\end{lemmaP}
\begin{proof}
 Let $g$ be the secondary Gauss map.
 Since the end is parabolic, the Schwarzian derivative $P:=S(g)$
 determines a $\PSU_{1,1}$-projective connection of parabolic monodromy.
 Then by \ref{item:parabolic} in Proposition~\ref{N} for 
 $\varepsilon=-1$ and $t=\pi$,
 there exists a $\PSU_{1,1}$ lift $g_0$ such that
 \[
     h(z) = \hat g_0(z) + i \log z\qquad
     \left( \hat g_0(z)=R\star g_0(z)=
            \frac{1}{i}
            \frac{g_0(z)+1}{g_0(z)-1}\right)
 \]
 is a meromorphic function on $\Delta^*$.
 Here, there exists a matrix $A\in\SU_{1,1}$ such that
 $g = A\star g_0$ or $1/g=A\star g_0$ holds.
 Thus, by the $\SU_{1,1}$-ambiguity of the secondary Gauss map,
 we may assume $g=g_0$ or $1/g_0$.
 Moreover, replacing $f$ with $-f$ if necessary
 (see Remark~\ref{rem:action-null}),
 we may assume $g=g_0$ without loss of generality.

 Since the end is $g$-regular, \ref{item:parabolic:2} 
 of Proposition~\ref{N} yields that $h(z)$ is meromorphic at $z=0$.
 Thus, we can write
 \[
     \hat g(z) = \hat g_0(z)
            = -i\log z + z^m\varphi(z)\qquad
           (\varphi(0)\neq 0, m\in\Z),
 \]
 where $\varphi(z)$ is a holomorphic function on a neighborhood
 of the origin.
 Then there exist an $a \in \R \setminus \{ 0 \}$ and
 a $\gamma\in (-\pi,\pi)$ such that
 \[
     \Im \hat g(z)=
         - \log r+a {r^m}\sin\bigl( m \theta+\gamma \bigr)
         + o(r^{m+1})  , \quad z=r e^{i \theta}.
 \]
 Here, the singular set $\{|g|=1\}$ is written as $\{\Im\hat g=0\}$.

 If the end is of the first kind, then $m\ge 0$ by 
 \ref{item:parabolic:3} in Proposition~\ref{N}. 
 Therefore, for each fixed $\theta$,
 the right-hand side approaches $\infty$ as $r\to 0$, which implies
 that the singular points do not accumulate at the end.

 If the end is of the second kind, then $m<0$.
 Therefore, for each fixed $\theta$,
 the right hand side approaches $\infty$ if $a\sin(m \theta+\gamma)>0$
 and $-\infty$ if $a\sin(m \theta+\gamma)<0$ as $r \to 0$,
 giving solutions of $\Im \hat g(z)=0$ for
 sufficiently small $r$ near the lines $\sin(m \theta+\gamma)=0$.
 This implies the second assertion.
\end{proof}
\begin{lemmaH}
 Let $f:\Delta^*\to S^3_1$ be a $g$-regular hyperbolic end.
 Then any ray in $\Delta^*$ emanating from the origin meets the
 singular set infinitely many times.
 {\rm (}See Figure~$\ref{fg:sing-h}$.{\rm )}
\end{lemmaH}
\begin{remark}
 This intersection property does not depend on the choice of a complex
 coordinate
 for a punctured neighborhood of the end.
\end{remark}
\begin{proof}[Proof of Lemma H]
 By Proposition~\ref{N} and an appropriate choice of $g$, 
 if we set $\hat g =R\star g$,
 then there is a $\mu \in \R \setminus \{ 0 \}$ such that
 \[
       h(z):=z^{-i\mu}\hat g(z)
 \]
 is a meromorphic function on $\Delta^*$. 
 Since $f$ is $g$-regular, 
 \ref{item:hyperbolic:2} of Proposition \ref{N} 
 implies $h$ has at most a pole at the origin, and
 $S(g)$ has a pole of order exactly $2$ at $z=0$.
 Thus we can rewrite
 \[
       \hat g(z)=z^{m+i\mu}\phi(z)
            \qquad (\phi(0)\ne 0, m \in \Z),
 \]
 where $\phi(z)$ is a single-valued holomorphic function on  
 $\Delta=\Delta^*\cup\{0\}$.
 Now, we set
 \[
       w=z \exp\left(\frac{\log\varphi}{m+i\mu}\right),
 \]
 which gives a new coordinate $w$ around the end, now at $w=0$.
 Then $\hat g(w)=w^{m+i\mu}$.
 Since $g=(\hat g-i)/(\hat g+i)$, 
 setting $w=r e^{i \theta}$, the singular set is
 \begin{align*}
  \{w\,;\, |g(w)|=1\}
  &=
  \{w\,;\, \Im(\hat g(w))=0\} \\
  &=
  \left\{(r,\theta)\in (0,1)\times (-\pi,\pi) \,;\,
     \mu \log r+ m\theta \equiv 0 \pmod{\pi}
  \right\},
 \end{align*}
 that is, 
 $r=\exp \bigl((n\pi-m\theta)/\mu\bigr)$, $n \in \Z$,
 which is a log-spiral when $m\ne 0$.
 If $m=0$, the singular set is a union of
 infinitely many disjoint circles.
 In any case, the singular set meets any ray based at $w=0$ infinitely 
 many times.  
 (See Figure~$\ref{fg:sing-h}$, left-hand side for the case $m=0$ and 
       right-hand side for the case $m\ne 0$.)
\end{proof}
\begin{corollary}\label{h-monodromy}
 The monodromy of a hyperbolic metric on $\Delta^*$  is either elliptic or
 parabolic.
 That is,
 hyperbolic monodromy never occurs.
\end{corollary}
\begin{proof}
 Suppose that a hyperbolic metric $d\sigma^2$ on $\Delta^*$ has
 hyperbolic monodromy.
 Let $g$ be a developing map for $d\sigma^2$.
 The data $(g,\omega=dz)$ produces an $F$ as in \eqref{eq:ode}, and then
 the immersion $f=F e_3 F^*: \Delta^* \to S_1^3$ is without singularities,
 since $d\sigma^2$ is nonsingular.  
 By Theorem \ref{thm:key}, $f$ is
 $g$-regular at $z=0$.  
 Then by Lemma H, the singularities accumulate at the end, a contradiction.
\end{proof}
Lemmas E1, E2, P and H imply:

\begin{corollary}[Characterization of hyperbolic ends]
 A $g$-regular end $f:\Delta^*\to S^3_1$
 of a  \cmcone{} face is hyperbolic if and only if
 every ray in $\Delta^*$ emanating from the origin meets the singular set
 infinitely many times.
\end{corollary}

\subsection*{Completeness}
We now give two theorems on complete \cmcone{} faces.
\begin{theorem}\label{thm:completeend}
 Any complete end of a \cmcone{} face is
 either $g$-regular elliptic or $g$-regular parabolic of the first kind.
\end{theorem}
\begin{proof}
 By Proposition~\ref{prop:complete}, the end is of puncture-type.
 Moreover, Theorem~\ref{thm:key} implies the end is
 $g$-regular. Thus the theorem follows from Corollary~\ref{h-monodromy}
 and Lemma P.
\end{proof}

\begin{theorem}\label{Thm:labeledbyYang3}
 Any complete \cmcone{} face is of finite type.  
 $($The definition of finite type is 
 given in Definition~\ref{def:ft}.$)$
\end{theorem}

\begin{proof}
 Let $f:M^2\to S^3_1$ be a complete \cmcone{} face.
 Then by Proposition~\ref{prop:complete}, 
 there is a compact Riemann surface
 $\overline M^2$ such that $M$ is biholomorphic to $\overline
 M^2\setminus\{p_1,\dots,p_n\}$.

 We fix any end $p_j$, and take a small coordinate
 neighborhood $(U,z)$ with $z=0$ at $p_j$.
 We may assume that
 there are no singular points on $U\setminus \{p_j\}$,
 and thus we may also assume that $|g|<1$ 
 on  $U\setminus \{p_j\}$ for a secondary Gauss map $g$.
 We know from Theorem~\ref{thm:completeend} that the end is
 a $g$-regular elliptic end or a
 $g$-regular parabolic end of the first kind.

 First, we consider the elliptic case.
 Since $|g|<1$, 
 we may assume that
 there exist some $\mu(\ge 0)$ 
 and a holomorphic function $h(z)$ on $U$
 with $h(0) \neq 0$ such that
 $g(z)=z^\mu h(z)$.
 If $|g(0)|=1$, then
 $\mu=0$ and 
 \[
    g(z)= e^{i\theta}\bigl(1+az^m+ o(z^{m})\bigr)
 \]
 for some $\theta \in \R$, $a\in \C\setminus\{0\}$ and
 $m\in \Z_+$, 
 which contradicts the fact that $|g|<1$
 on $U \setminus \{ p_j \}$.
 Hence $|g(0)| < 1$. 
 Therefore, there exist a neighborhood $\widetilde{U} \subset U$ of
 $p_j$ and
 an $\varepsilon > 0$
 such that $|g|^2<1-\varepsilon$ on $\widetilde{U}$. 
 So, on $\widetilde{U}$,
 \begin{equation}\label{eq:c-estimate}
    K_{ds^2}\,ds^2 = \frac{4 |dg|^2}{(1-|g|^2)^2}
           \leq \frac{4}{\varepsilon^2}
            \frac{4\,|dg|^2}{(1+|g|^2)^2}
           = \frac{4}{\varepsilon^2} (-K_{d\hat s^2})\,d\hat s^2,
 \end{equation}
 where $d\hat s^2$ is the metric as in \eqref{eq:dshat^2},
 which is defined on $\widetilde{U}\setminus\{p_j\}$
 because $g(z)=z^{\mu}h(z)$.
 Since $p_j$ is a regular singularity
 of the punctured spherical metric
 \[
    d\hat\sigma^2 := (-K_{d\hat s^2})\,d\hat s^2
          =\frac{4\,|dg|^2}{(1+|g|^2)^2},
 \]
 $d\hat\sigma^2$ has finite area, so $ds^2$ has finite total curvature on 
 $\widetilde{U}\setminus\{p_j\}$, by \eqref{eq:c-estimate}.

 Next we consider the parabolic case.
 By Theorem~\ref{thm:completeend}, 
 the end is parabolic of the first kind.
 Then  without loss of generality, we may assume
 there exists a holomorphic
 function $h(z)$ on $\Delta$ such that
 (we set $t=\pi$ and replace $h(z)$ by $ih(z)$
 in the proof of Proposition~\ref{N})
 \begin{equation}\label{eq:para-first}
    \hat g(z)=i(h(z) \pm \log z), \qquad\text{where}\quad
      \hat g(z)=R\star g(z)=
       \frac{1}{i}\frac{g(z)+1}{g(z)-1}.
 \end{equation}
 If we set $k(z):=h(z) \pm \log z+1$, we have
 $g = 1-2/k$, $g' = 2k'/k^2$ and 
 \[
  1-|g(z)|^2=\frac{4(\Re k(z)-1)}{|k(z)|^2}
            =\frac{4(\Re h(z) \pm \Re\log z)}{|k(z)|^2}.
 \]
 So we have
 \begin{equation}\label{eq:p-metric}
  d\sigma^2=\frac{4\,|dg|^2}{(1-|g|^2)^2}
   =\frac{|h'(z) \pm (1/z)|^2|dz|^2}{(\Re h(z) \pm \Re\log z)^2}.
 \end{equation}
 We set $c:=\sup_{z\in \Delta}|\Re h(z)|$ and $r = |z|$.
 Since $\log r\to -\infty$ as $z \to 0$, we may assume $-\log r > c$.  Then
 \[
    |\Re h(z) \pm \log r|\ge
        \bigl||\Re h(z)|-|\log r|\bigr|\ge |c+\log r|, \quad \text{and}
 \]
 \begin{equation}\label{eq:parab-met-estimate}
      d\sigma^2\le \frac{C^2}{r^2 (c+\log r)^2} |dz|^2,
 \end{equation}
 where $C=\sup_{z\in \Delta}|zh'(z)+1|$.  Since 
 \[
    \int_0^\varepsilon \frac{C^2 r\,dr}{r^2(c+\log r)^2 }
      =-\frac{C^2}{c+\log \varepsilon}<\infty,
 \]
 the area of a sufficiently small
 punctured neighborhood of $z=0$ with respect to $d\sigma^2$ is finite,
 which proves the assertion.
\end{proof}
Theorem~\ref{thm:main} in the introduction 
follows from Theorems~\ref{thm:completeend} and
\ref{Thm:labeledbyYang3}.

\section{The light-cone Gauss map 
and\\extrinsic behavior of ends}
\label{sec:lgm}
Let $\LC=\{x\in \R^4_1\,;\,\inner{x}{x}=0\}$ be the light-cone of
$\R^4_1$, with future and past light cones
\[
   \LC_{\pm} :=\{x=(x_0,x_1,x_2,x_3)\in \LC\,;\,\pm x_0>0\}.
\]
The multiplicative group $\R_+$ of the positive real numbers
acts on $\LC_{\pm}$ by scalar multiplication.
The {\em ideal boundary\/} $\partial S^3_1$ of $S^3_1$ consists of two
(future and past) components
\[
   \partial_\pm S^3_1:= \LC_{\pm}/\R_+,
\]
each of which are identified with $\C\cup\{\infty\}$ by the projection
\begin{equation}\label{eq:proj-lc}
   \pi \colon \partial_{\pm} S^3_1
   \ni \bigl[(v_0,v_1,v_2,v_3)\bigr]\longmapsto
       \frac{1}{v_0-v_3}(v_1+i v_2)\in \C\cup\{\infty\}.
\end{equation}
The isometries of $\R^4_1$
induce M\"obius transformations on $\C\cup\{\infty\}$.
The boundary $\partial S^3_1$ is identified with the set of 
equivalence classes of oriented time-like geodesics in $S^3_1$.

In particular, for a space-like immersion $f\colon{}M^2\to S^3_1$
with the (time-like) unit normal vector field $\nu$,
the equivalence class $[f+\nu]$ determines a point in $\partial S^3_1$
for each $p\in M^2$.
Hence we have the {\em light-cone Gauss map}
\[
    L=[f+\nu]\colon{}M^2 \longrightarrow \partial S^3_1.
\]

Let $f\colon{}M^2\to S^3_1$ be a \cmcone{} face, and $p\in M^2$
a regular point, that is, $f$ is an immersion in a neighborhood
of $p$.  
Under the identification of $\R^4_1$ and $\Herm(2)$ as in 
\eqref{eq:herm-mink}, we can compute that
the unit normal vector $\nu$
is
\[
    \nu = \frac{1}{|g|^2-1} F
    \begin{pmatrix}
        1+|g|^2  & 2 g \\
        2 \bar g & 1 + |g|^2
    \end{pmatrix} F^*,
\]
where $F$ is the holomorphic lift of $f$ and 
$g$ is the secondary Gauss map.
Hence
\begin{equation}\label{eq:lgm}
   L = \!\left[
              \frac{2}{|g|^2-1}
          F
          \begin{pmatrix}
           |g|^2  &  g \\
            \bar g & 1
          \end{pmatrix} F^*
             \right]
           = \sign (|g|^2-1)\!\!
             \left[
          F
          \begin{pmatrix}
           |g|^2  &  g \\
            \bar g & 1
          \end{pmatrix} F^*
             \right],
\end{equation}
where $\sign (|g|^2-1)$ is the sign of the function
$|g|^2-1$.  Thus, we have:
\begin{proposition}\label{prop:lgm}
 The light-cone Gauss map $L$ of a \cmcone{} face takes values in
 $\partial_+S^3_1$ {\rm (}resp.~$\partial_- S^3_1${\rm )}
 if $|g|>1$ {\rm (}resp. $|g|<1${\rm )}.
 Moreover, its projection $\pi\circ L$ is the
 hyperbolic Gauss map $G$ as in \eqref{eq:h-gauss},
 which extends to the singular set.
\end{proposition}
\begin{proof}
 By \eqref{eq:lgm},
 the $x_0$-component of $f+\nu$ is
 \[
 \frac{1}{|g|^2-1}
 \trace  \left(F
 \begin{pmatrix}
           |g|^2  &  g \\
            \bar g & 1
 \end{pmatrix} F^* \right)
  =
 \frac{1}{|g|^2-1}\left(|gF_{11}+F_{12}|^2+|gF_{21}+F_{22}|^2\right),
 \]
 where $F=(F_{ij})_{i,j=1,2}$.
 Here,
 $|gF_{11}+F_{12}|^2+|gF_{21}+F_{22}|^2>0$ holds because
 $\det F=1\neq 0$, implying the first part of the proposition.
 By \eqref{eq:lgm}, \eqref{eq:proj-lc},
 \eqref{eq:ode} and \eqref{eq:h-gauss}, we have
 \[
    \pi\circ L =  \pi\circ[f+\nu]= \frac{gF_{11}+F_{12}}{gF_{21}+F_{22}}
                   = \frac{dF_{11}}{dF_{21}}=G, 
 \]
 and this completes the proof.  
\end{proof}
Next we give a criterion for when a complete regular end approaches
$\partial_+S^3_1$ or $\partial_-S^3_1$:

\begin{proposition}\label{prop:fatend}
 Let
 $f:\Delta^*\to S^3_1$ be a complete regular end at $z=0$
 and let $g$ be a secondary Gauss map of $f$.
 Then the image of $f$ converges to a point in $\partial_+S^3_1$
 {\rm (}resp. $\partial_-S^3_1${\rm )}
 at the end if and only if $|g|<1$ {\rm (}resp. $|g|>1${\rm )} near the end.
\end{proposition}
\begin{proof}
 We can change the holomorphic null
 lift $F$ to $F^{\natural}$ as in
 \eqref{eq:time-reverse-lift}, so that $f$ and its
 secondary Gauss map $g$ change to $-f$ and
 $1/g$.  
 The end of $f$ approaches $\partial_{\pm}S^3_1$ if and only if
 $-f$ approaches $\partial_{\mp}S^3_1$, so
 it is sufficient to prove this result under
 the assumption $|g|<1$ on $\Delta^*$.
 By Theorem~\ref{thm:completeend}, the end is either elliptic or parabolic.

 First we assume the end is elliptic.
 Replacing $F$ by $aFb^{-1}$ ($a\in\SL_2\C$, $b\in\SU_{1,1}$) and
 using the Weierstrass preparation theorem if necessary,
 we may assume without loss of generality that the hyperbolic and secondary
 Gauss maps are
 \[
    G(z) =z^m, \qquad g(z)=z^\mu h(z)
    \qquad \text{for some } m\in \Z_+,~\mu\in\R\setminus\{0\} \; , 
 \]
 where $h$ is a holomorphic function on 
 $\Delta$ with $h(0)\ne 0$.  Here $\mu>0$ because $|g|<1$.

 If $m\neq \mu$, Small's formula \eqref{eq:small} implies that
 \[
     F = \frac{1}{2\sqrt{m\mu}}
      \begin{pmatrix}
       -z^{\frac{\hphantom{-}m-\mu}{2}}(m+\mu) \bigl(1+o(1)\bigr) &
    \hphantom{-} z^{\frac{\hphantom{-}m+\mu}{2}}(m-\mu) \bigl(1+o(1)\bigr) \\
       \hphantom{-} z^{\frac{-m-\mu}{2}}(m-\mu) \bigl(1+o(1)\bigr) &
       -z^{\frac{-m+\mu}{2}}(m+\mu) \bigl(1+o(1)\bigr)
      \end{pmatrix}.
 \]
 Since $m$ and $\mu$ are positive, the first component $x_0$ 
 is 
 \[
    x_0 = \frac{1}{2}\trace (Fe_3 F^*)=
           \frac{(m-\mu)^2}{8m\mu}r^{-m-\mu}\bigl(1+o(1)\bigr)
            \to +\infty\qquad (r\to 0),
 \]
 where $z=re^{i\theta}$.
 Other components of $f=(x_0,x_1,x_2,x_3)$ are expressed as 
 \begin{align*}
    x_1 + i x_2 &= 
     e^{im\theta}r^{-\mu}\frac{\mu^2-m^2}{4\mu m}\bigl(1+o(1)\bigr),\\
    x_3 &= -\frac{(m-\mu)^2}{8m\mu}r^{-m-\mu}\bigl(1+o(1)\bigr).
 \end{align*}
 We now consider the stereographic projection given in \cite{F}:
\begin{multline}\label{eq:proj}
   \Pi:\{(x_0,x_1,x_2,x_3)\in S^3_1\,;\, x_0>1\}
       \ni(x_0,x_1,x_2,x_3)\\
   \longmapsto\frac{1}{1+x_0} (x_1,x_2,x_3)\in 
       \left\{(X_1,X_2,X_3)\in \R^3\,;\,
       \frac{1}{2}<\sum_{j=1}^3(X_j)^2<1 \right\} ,
\end{multline}
  which is a diffeomorphism.
 Then $\Pi\circ f$ is expressed as
 \[
  \Pi\circ f = (0,0,-1) + o(1).
 \]
  Thus, $\Pi\circ f$ approaches $(0,0,-1)\in S^2=\partial_+ S_1^3$.

 When $\mu=m$, by \eqref{eq:small} again,
 $F_{11}$, $F_{12}$ and $F_{22}$ are bounded
 on a neighborhood of $0$, and 
 these components can be extended to 
 become holomorphic
 on a neighborhood of $0$.
 If $F_{21}$ is bounded, $F$ must be holomorphic and then
 the induced metric is bounded, which contradicts the 
 weak completeness  of the end.
 Hence $F_{21}$ has a pole at $0$.
 So the $x_0$-component of $f$ is 
 \[
   x_0 = \frac{1}{2}|F_{21}|^2 + \text{(a bounded function)}
       \to +\infty\qquad (z\to 0).
 \]
 Moreover, since
 \begin{align*}
     x_3 &= -\frac{1}{2}|F_{21}|^2 + \text{(a bounded function)},\\
     x_1 + i x_2 &= F_{11}\overline{F_{21}}-F_{12}\overline{F_{22}}
         = c \overline{F_{21}} + \text{(a bounded function)},
 \end{align*}
 we have $\Pi\circ f \to (0,0,-1)$ as $z\to 0$.

 Next we assume the end is parabolic.
 Again we may set $G=z^m$, $m\in \Z_+$.
 Applying 
 Proposition~\ref{N} for $t=2m\pi$ and $\varepsilon=-1$,
 there exists a $\PSU_{1,1}$-lift $g_0$
 of $S(g)$ such that 
 $h(z):=R\star g_0(z)+2mi\log z$ is a single-valued meromorphic function
 on $\Delta^*$, and the secondary Gauss map $g$ satisfies
 $g=A\star g_0$ or $1/g=A\star g_0$ for some $A\in\SU_{1,1}$.
 By completeness, Lemma~P  implies that the end is of first kind.
 Hence by \ref{item:parabolic:3} of Proposition~\ref{N},
 $h(z)$ is holomorphic on a neighborhood of the origin.
 Moreover, by the assumption $|g|<1$, \ref{item:parabolic:4}
 of Proposition~\ref{N} yields that $g=A\star g_0$ 
 for some $A\in\SU_{1,1}$.
 Thus, without loss of generality, we may set
 \[
    R\star g(z) =2mi\bigl(k(z)-\log z\bigr),
 \]
 here we set $h(z)=2mik(z)$.

 For a holomorphic null lift $F$ of $f$ with the secondary Gauss
 map $g$, set
 \[
     \hat F = FB^{-1},\qquad \text{where}\quad
     B: = 
        \left( \begin{pmatrix} 1 & 0 \\ 0 & 2i\end{pmatrix} R \right) 
       =\begin{pmatrix}1/2 & 1/2\\ -1&1\end{pmatrix}
        \in\SL_2\C.
  \]
 Then $\hat F$ is a holomorphic null immersion 
 whose
 hyperbolic Gauss map $\hat G$  and secondary Gauss map $\hat g$
 are given by
 \begin{equation}\label{eq:parabolic-modified-Gg}
  \begin{aligned}
    \hat G &= G = z^m,\\
    \hat g(z) &= B\star g(z)  
      = \begin{pmatrix} 1 & 0 \\ 0 & 2 i \end{pmatrix}\star 
          \bigl(R\star g(z)\bigr) = 
           m\bigl(k(z)-\log z \bigr).
  \end{aligned}
 \end{equation}
 So applying \eqref{eq:small} for this $(\hat G,\hat g)$,
 the components of $\hat F$ are written as 
 \begin{equation}\label{eq:F-hat}
  \begin{alignedat}{2}
   \hat F_{11}(z) &= -\frac{i}{2} z^{m/2}\varphi_1(z),\qquad
   &\hat F_{12}(z) &= -\frac{i}{2} z^{m/2}
         \bigl(m\varphi_1(z)\log z+\psi_1(z)\bigr),\\
   \hat F_{21}(z) &= 
            \frac{i}{2}
            z^{-m/2}\varphi_2(z),\qquad
   &\hat F_{22}(z) &= 
            \frac{i}{2}z^{-m/2}\bigl(m\varphi_2(z)\log z+\psi_2(z)\bigr),
  \end{alignedat}
 \end{equation}
 where $\varphi_1$, $\varphi_2$, $\psi_1$ and $\psi_2$ are holomorphic
 functions defined on a neighborhood of the origin 
 such that
 \[
     \varphi_1 (0) = \varphi_2 (0) =1 .
 \]
 Since $Be_3B^*=-e_1$, $f=Fe_3F^*$ satisfies
 \begin{equation}\label{eq:parabolic-f-repr}
    f= -\hat F \begin{pmatrix}0 & 1 \\ 1 & 0\end{pmatrix} \hat F^*
      = -\begin{pmatrix}
      \hat F_{11}\overline{\hat F_{12}} +
      \hat F_{12}\overline{\hat F_{11}}     &
      \hat F_{11}\overline{\hat F_{22}} +
      \hat F_{12}\overline{\hat F_{21}} \\
      \overline{\hat F_{11}}{\hat F_{22}} +
      \overline{\hat F_{12}}{\hat F_{21}} &
      \hat F_{21}\overline{\hat F_{22}} +
      \hat F_{22}\overline{\hat F_{21}}
    \end{pmatrix}.
 \end{equation}
 Hence the components of $f$ are expressed as
 \begin{equation}\label{eq:parabolic-end}
  \begin{aligned}
   x_0 &= \hphantom{-}\frac{m}{4} r^{-m}\bigl(\hphantom{-}
           \eta_1(u,v)\log r+\delta_1(u,v)\bigr),\\
   x_3 &= \hphantom{-}\frac{m}{4} r^{-m}
    \bigl(-\eta_2(u,v)\log r+\delta_2(u,v)\bigr),\\
   x_1 + i x_2 &=
      -\frac{m}{2}e^{im\theta}\bigl(
               \hphantom{-}\eta_3(u,v)  \log r + \delta_3(u,v)
                              \bigr),
  \end{aligned}
 \end{equation}
 where $z=re^{i\theta}=u+iv$.
 Here,  $\eta_j(u,v)$ ($j=1,2$) and $\delta_j(u,v)$ ($j=1,2$) 
 (resp.\ $\eta_3(u,v)$ and $\delta_3(u,v)$)
 are real-valued (resp.\ complex-valued)  differentiable functions 
 defined on a neighborhood of the origin, such that
 $\eta_j(0,0)=1$ ($j=1,2,3$).

 The equations \eqref{eq:parabolic-end} yield that
 $x_0\to +\infty$ and  $\Pi\circ f\to(0,0,-1)$ as $z\to 0$.
\end{proof}

\section{The Osserman type inequality}
\label{sec:fujimori}

Here we prove Theorem~\ref{thm:fujimori}
stated in the Introduction.  First we prepare:

\begin{lemma}\label{lem:complete-p}
 The Hopf differential of a \cmcone{} face
 has a pole of order $2$ at any complete regular parabolic end.
\end{lemma}
\begin{proof}
 Let $f:\Delta^*\to S^3_1$ be a complete regular parabolic end at $z=0$. 
 By Theorem \ref{thm:completeend}, the end is of the first kind.
 Then
 \[
   2Q+S_z(G)\,dz^2 =S_z(g)\,dz^2 =\frac{1}{z^2}
       \left( \frac{1}{2} +o(1)\right)\,dz^2.
 \]
 Since $G$ is meromorphic at $z=0$, 
 we may assume that $G=z^m\varphi(z)$, where $m$ is a positive integer
 and $\varphi(z)$ is a holomorphic function on a neighborhood 
 of $0$ with $\varphi(0)\neq 0$.  
 Applying \eqref{eq:s_ord} to $S_z(G)$, it follows that
 $Q$ has a pole of order $2$ at $z=0$.
\end{proof}

The next lemma improves a result in \cite[Proposition 4.4]{F}:

\begin{lemma}\label{lem:osserman-ineq}
 Let $f:\Delta^*\to S^3_1$ be a complete regular end at $z=0$ of a
 \cmcone{} face
 with Hopf differential $Q$ and hyperbolic Gauss map $G$.
 Then the ramification order $m$ of $G(z)$ at $z=0$ 
 satisfies
 \begin{equation}\label{eq:local-oss}
  m\ge \Ord_{z=0}(Q)+3 \; .
 \end{equation}
 {\rm(}For the definition of the ramification order, see
 the subsection about the Schwarzian derivative in
 Section~\ref{sec:prelim}.{\rm )}
 Here, $\Ord_{z=0}Q$ denotes the order of $Q$ at the origin, 
 that is, $\Ord_{z=0}Q=k$ if $Q=z^k\varphi(z)\,dz^2$, where $\varphi(z)$ is 
 holomorphic at $z=0$ and $\varphi(0)\neq 0$.
\end{lemma}
\begin{proof}
 By Theorem~\ref{thm:completeend}, a complete end is either elliptic or
 parabolic of the first kind.
 The elliptic case has been proved in \cite{F}.
 Assume that the end is parabolic. 
 Then by Lemma \ref{lem:complete-p},
 $Q$ must have a pole of order $2$ at $z=0$,
 which proves the inequality since $m\geq 1$.
\end{proof}
It should be remarked that the order of the metric 
$d\sigma^2_{\#}=4|dG|^2/(1+|G|^2)^2$ at $0$ 
is equal to $m-1$, 
where $m$ is the ramification order of $G$.
Using Lemma \ref{lem:osserman-ineq} instead of \cite[Proposition 4.4]{F}, 
the inequality in Theorem~\ref{thm:fujimori}
is proved in the same way as \cite{UY5}, \cite{F}.

The condition for equality in \eqref{eq:fujimori} 
in Theorem~\ref{thm:fujimori}
for elliptic ends was completely analyzed in \cite{F}.
So, it suffices to show the following theorem for 
parabolic ends.
Note that $\Ord_p(Q)=-2$ for complete regular parabolic ends,
hence the equality in \eqref{eq:local-oss} holds if and only if 
$G$ does not branch at $p$
(see \cite{F} for details). 
\begin{theorem}\label{thm:equality}
 A complete regular parabolic end of a  \cmcone{} face is 
 properly embedded if and only if the hyperbolic Gauss map $G$ 
 does not branch at the end.
\end{theorem}
\begin{proof}
 Let $f:\Delta^*\to S^3_1$ be a complete regular parabolic
 end at $z=0$.
 Taking $-f$ instead of $f$ if necessary, we may assume that $|g|<1$ in
 a neighborhood of the end, and that 
 $G(z) = z^m, m \ge 1$ and 
 $g(z) = R^{-1} \star \left( 2mi (k(z) - \log z) \right)$,
 as in the proof of Proposition~\ref{prop:fatend}.
 Then $f$ is represented as in 
 \eqref{eq:F-hat} and \eqref{eq:parabolic-f-repr}.

 By Proposition~\ref{prop:fatend}, the image of $f$ tends to
 a point in
 $\partial_+ S^3_1$.
 So we may assume that $x_0>1$ on $\Delta^*$, and 
 \[
    \Pi\circ f\colon{}\Delta^*\ni z\longmapsto (X_1,X_2,X_3)\in\R^3
 \]
 is well-defined, 
 where $\Pi$ is the projection in \eqref{eq:proj}.

 Here, by \eqref{eq:parabolic-end}, 
 \[
   U(z) := z^{-m}\frac{x_1+i x_2}{1+x_0}
    \left(\vphantom{\int}=z^{-m}(X_1+iX_2)\right)
       =
    -2 \frac{\eta_3 \log r+\delta_3}{\eta_1\log r + \delta_1+(4/m)r^m}.
 \]
 Since $\eta_1,\eta_3$ and $\delta_1,\delta_3$ are 
differentiable functions
 defined on a neighborhood of $0$, 
 we have
 \begin{equation}\label{eq:U}
  \lim_{z\to 0}U(z)=-2\ne 0, \;\; 
   \lim_{z\to 0}z \frac{\partial}{\partial z} U(z)=0
  \;\; \text{and}\;\;
  \lim_{z\to 0}z \frac{\partial}{\partial \bar z} U(z)=0.
 \end{equation}

 Now we suppose that the ramification order $m$
 of the hyperbolic Gauss map at $z=0$ is $1$, 
 that is $m=1$.
 As seen in the proof of Proposition~\ref{prop:fatend}, $\Pi\circ f$ 
 converges to $(0,0,-1)$.
 Then $X_1+i X_2=zU(z)$ and \eqref{eq:U} yield that
 \[
   \lim_{z\to 0}\frac{\partial}{\partial z} (X_1+i X_2) \ne 0
   \qquad \text{and}\qquad
   \lim_{z\to 0}\frac{\partial}{\partial \bar z} (X_1+i X_2) = 0 \; ,
 \]
 which implies that the correspondence
 $z\mapsto X_1+iX_2$
 is bijective near the origin, and the end is properly embedded.

 Conversely, suppose that the end is properly embedded.
 We have already seen that $X_3\to -1$ as $z \to 0$.
 Moreover $U(0)\ne 0$ implies that for any
 sufficiently small $\varepsilon>0$, the image of
 the end $f(\{z ; 0 < |z| < \epsilon \})$ does not
 meet the $X_3$-axis and is diffeomorphic to a cylinder.
 Then the image of $\Pi \circ f(\{z ; |z| = \epsilon \})$ 
 by the orthogonal projection
 $(X_1,X_2,X_3) \mapsto X_1+iX_2$ is an embedded closed curve
 with the winding number $m$ with respect to the origin.
 So $m=1$.
\end{proof}
\begin{remark}
 In Proposition 4.4 of \cite{F}, the first 
 author showed the equality condition 
 in Theorem~\ref{thm:fujimori} for elliptic ends
 using the expression of the solution of
 the ordinary differential equation \eqref{eq:ode}.  
 Here we proved the equality condition in Theorem~\ref{thm:fujimori} for
 parabolic 
 ends by using Small's formula \eqref{eq:small}.  It is also 
possible to 
 prove the result in \cite{F} 
 more directly by using  \eqref{eq:small}. 
\end{remark}

We give here three important examples:
\begin{example}[%
An incomplete 3-noid not satisfying \eqref{eq:fujimori} in
Theorem \ref{thm:fujimori}]
 We set
 $M^2 :=\C\setminus\{0,1\}$
 and
 \[
   G:=z,\qquad
   g:=\frac{2z-1}{2z(z-1)}-\log\frac{z}{z-1} \; .
 \]
 Then \eqref{eq:small} gives a \cmcone{} face
 $f:M^2\to S^3_1$
 with hyperbolic and secondary Gauss maps $G$ and $g$, and Hopf differential
 \[
    Q=\frac{1}{2}(S(g)-S(G))= -\frac{2 dz^2}{z(z-1)}.
 \]
 Since the lift metric
 \[
   ds^2_\# = \frac{4(1+|z|^2)^2}{|z(z-1)|^2}|dz|^2
 \]
 is complete on $M^2$, $f$ is weakly complete.
 The end $z=\infty$ is complete and elliptic,
 and $z=0,1$ are parabolic ends of the second kind.
 Hence $z=0,1$ are incomplete ends.
 Since $\deg(G)=1$, $f$ does not satisfy \eqref{eq:fujimori}.
 This implies that completeness is an essential
 assumption in Theorem~\ref{thm:fujimori} in the introduction.
\end{example}
\begin{example}[%
A 2-noid with complete parabolic ends 
 satisfying the equality in \eqref{eq:fujimori}]
We set
\begin{equation}
 F(z)=
   \dfrac{i}{2\sqrt{2}}
   \begin{pmatrix}
    \sqrt{z} & 0 \\ 0 & \sqrt{z}^{-1}
   \end{pmatrix}\!
  \begin{pmatrix}
    3-\log z & -1+\log z \\
   1+\log z & -3-\log z
  \end{pmatrix}.
\end{equation}
 Then $f=Fe_3F^*:\C\setminus \{0\}\to S^3_1$ has
 two parabolic regular ends.
 The hyperbolic Gauss map $G$, the secondary Gauss map $g$ 
 and the Hopf differential $Q$ are computed as follows: 
\[
 G=z, \qquad  
 g=\frac{\log z+1}{\log z-1}, \qquad
 Q=\frac{dz^2}{4z^2}.  
\]
 Since $\{z\in\C\,;\,|g(z)|=1\}=\{z\in\C\,;\,|z|=1\}$, 
 the singular set is compact, and hence $f$ is complete. 

 Any genus zero \cmcone{} face with two parabolic regular ends 
 and with degree $1$ hyperbolic Gauss map is congruent to this $f$.
 We call this \cmcone{} face the {\em parabolic catenoid}.
 On the other hand, the \cmcone{} face with 
 $G=z$, $g=z^{\mu}$ ($\mu\in\R\setminus\{0\}$)
 given in \cite[Example 5.4]{F} is called the {\em elliptic catenoid}.
\end{example}

\begin{example} \label{ex:elliptic-integral}
 (A complete 4-noid with 4 integral elliptic ends
 satisfying the equality in \eqref{eq:fujimori})\par
 Since $\SL_2\C$ can be identified 
 with the complex hyperquadric $Q^3$ 
 of $\C^4$, the null (meromorphic) curves in $\SL_2\C$ can
 be identified with those in $Q^3$.
 The null curve in $\SL_2\C$ with
 \[
 G := \frac{3(z^3+2)}{4-z}, \qquad g:=-\frac{z^3-12z^2+2}{3z}
 \]
 belongs to the moduli space $\mathscr{M}_4$ in the classification 
 list of null curves in $Q^3$ in Bryant \cite{B2}, which
 has four integral elliptic ends at the roots of $1+6z^2-z^3$ and $z=\infty$.
 Since $G$ is of degree $3$ and $\chi(\C\cup \{\infty\})=2$,
 the corresponding \cmcone{} face attains equality in
 \eqref{eq:fujimori} of Theorem~\ref{thm:fujimori}.
 (For the definition of an integral elliptic end, see Definition 
 \ref{def:labeledbyYang2}.)
\end{example}

\begin{remark}
 We can deform an elliptic catenoid to a parabolic catenoid.
 Let $f_\mu$ be an elliptic catenoid with 
 the hyperbolic Gauss map $G=z$ and 
 the secondary
 Gauss map $g=z^\mu$, where $\mu>0$.
 Then the hyperbolic metric
 corresponding to $f_\mu$ is
 \[
    d\sigma_\mu^2=\frac{4 |dg|^2}{(1-|g|^2)^2}
      =\frac{4\mu^2|z|^{2\mu-2}}{(1-|z|^{2\mu})^2}|dz|^2 \; .
 \]
 It can be easily checked that
 \[
    \lim_{\mu\to 0} d\sigma^2_\mu
  =\frac{|dz|^2}{(r\log r)^2}, \quad \text{where}\quad z=re^{i\theta},
 \]
 which is the hyperbolic metric of a parabolic catenoid with
 \[
 g(z)= R \star \log z = \frac{1}{i}\frac{\log z+1}{\log z -1 },
 \]
 see \eqref{eq:p-metric}.
 On the other hand, by Small's formula \eqref{eq:small}
 there exists a unique smooth $1$-parameter family of
 \cmcone{} faces $\tilde f_\mu$ ($\mu\ge 0$)
 with hyperbolic Gauss map $G=z$
 and associated hyperbolic metric $d\sigma^2_\mu$.
 Then $\tilde f_\mu$ is congruent to $f_\mu$,
 and $\tilde f_0$ gives a parabolic catenoid.
\end{remark}

\begin{remark}
As a consequence of Remark \ref{rem:complete}, 
we know that there are no compact \cmcone{} immersed surfaces in
$S^3_1$.  Here we give an alternative proof of this: 
Let $M^2$ be a compact Riemann surface without boundary,
and suppose there exists a compact \cmcone{} face
$f:M^2\to S^3_1$ which has no singular points.
Let $F$ be a holomorphic lift of $f$.
We may assume that $|g|<1$ since there are no singular points.  
Then, by \eqref{eq:ode}, we have
\[
   f_{z\bar z}=(1-|g|^2)F 
   \begin{pmatrix}
    g \\ 1 
   \end{pmatrix}
   \begin{pmatrix}
    \bar g &  1 
   \end{pmatrix}
   F^*|\hat\omega|^2,
\]
where $z$ is a local complex coordinate and $\omega=\hat\omega\,dz$.
Thus, $\trace f$ is a nonconstant subharmonic function, 
which is a contradiction to the maximum principle.
\end{remark}

This proof does not apply to compact \cmcone{} faces,
leading us to the following open problem:
\begin{problem*}
 Is there a compact \cmcone{} face?
\end{problem*}
 
 If such a \cmcone{} face exists, the genus $\gamma$ 
 must be greater than or equal to $3$,
 since equality in \eqref{eq:fujimori} in the introduction
 holds in this case and
 the degree of the hyperbolic Gauss map must be $\gamma-1$.

\appendix
\section{Meromorphicity of the Hopf differential}
\label{app:Q-mero}
In this appendix, we shall give a proof of the following
\begin{theorem}\label{thm:Q-mero}
 Let $\overline M^2$ be a compact Riemann surface.
 Then the Hopf differential $Q$ of a complete \cmcone{} face
 \[
    f\colon{}\overline{M}^2\setminus\{p_1,\dots,p_n\}
       \longrightarrow S^3_1
 \]
 is meromorphic on $\overline{M}^2$.
\end{theorem}

\begin{proof}
 It is sufficient to show the meromorphicity of $Q$ 
 at a complete end
 $f\colon{}\Delta^*=\{z;0<|z|<1\}\to S^3_1$
 at the origin.
 We write the Hopf differential $Q$ as 
 \[
    Q = \hat Q \,dz^2,
 \]
 where $\hat Q$ is a holomorphic function on $\Delta^*$.
 By Theorem~\ref{thm:completeend},
 a complete end $f\colon{}\Delta^*\to S^3_1$ is either 
 a $g$-regular elliptic end or a 
 $g$-regular parabolic end of the first 
 kind. (The definition of $g$-regularity is given in Definition
 \ref{def:g-reg}.)

 First, we consider the case that $f$ is elliptic.
 By $g$-regularity, the secondary Gauss map is written in the form 
 \[
    g = z^{\mu} h(z)\qquad 
    (\text{$h$ is a holomorphic function with $h(0)\neq 0$}),
 \]
 where $\mu$ is a real number.
 Since $|g(0)|\neq 1$ by completeness, 
 we may set $g(0)=0$, or $\infty$, because of the 
 $\SU_{1,1}$-ambiguity of $g$.
 Moreover, replacing $f$ by $-f$ if necessary, we may assume
 $\mu>0$ without loss of generality.
 In this case, the corresponding hyperbolic metric $d\sigma^2$
 is written as
 \[
    d\sigma^2= \left(
                 \frac{2|z|^{\mu-1}|\mu h(z)+zh'(z)|}{%
                      \bigl|1-|z|^{2\mu}|h(z)|^2\bigr|}\,|dz|
               \right)^2,
               \qquad\left({~}'=\frac{d}{dz}\right).
 \]
 Since $|z|^{\mu}|h(z)|$ and $zh'(z)$ tend to $0$ as $z\to 0$
 and $h(z)$ is bounded near the origin, we have that 
 \[
     d\sigma \geq c |z|^{\mu-1}\,|dz| \geq c |z|^l \,|dz|
 \]
 holds on a neighborhood of the origin,
 where $l$ is the smallest integer such that $l\geq \mu-1$ 
 and $c$ is a positive constant.
 Then, by \eqref{eq:gauss}, we have
 \[
    ds = 2\frac{|Q|}{d\sigma}
       \leq  2 \frac{|\hat Q|}{c|z|^l}|dz|
       = \frac{2}{c}\left|
            \frac{\hat Q}{z^l}\,dz
         \right|.
 \]
 Since $ds$ is complete at $0$, we have meromorphicity 
 of the one-form $z^{-l}\hat Q \,dz$ at the 
 origin, because of \cite[Lemma 9.6, page 83]{O}.

 Next, we consider the case that $f$ is parabolic.
 Since the end is $g$-regular parabolic of the first kind,
 one can choose the secondary Gauss map $g$ as in \eqref{eq:para-first}:
 \[
   g=R^{-1}\hat g,\qquad \hat g(z)=i(h(z)\pm \log z),
 \]
 where $h(z)$ is a holomorphic function on
 $\Delta:=\Delta^* \cup \{ 0 \}$.
 Hence $d\sigma^2$ is written as in \eqref{eq:p-metric}:
 \[
   d\sigma^2 = \left(
                 \frac{|h'(z)\pm (1/z)|}{%
                   \bigl|\Re h(z)\pm \log |z|\bigr|
                 }\,|dz|
               \right)^2
               \qquad \left({~}'=\frac{d}{dz}\right).
 \]
 Since $z\log|z|\to 0$ as $z\to 0$ and $h$ is bounded
 on a neighborhood of the origin, 
 \[
   d\sigma = \frac{|1\pm zh'(z)|}{%
            \bigl|z\log|z|\bigr|
            \left|
               1\pm \frac{\Re h(z)}{\log |z|}
            \right|}\,|dz|
            \geq c|dz|
\]
holds on a 
neighborhood of the origin, where $c$ is a positive constant.
Thus, 
\[
    ds = 2\frac{|Q|}{d\sigma}
       \leq  \frac{2}{c}|\hat Q\, dz|.
\]
Hence, by the same argument as in the elliptic case,
we have meromorphicity of $Q$ at the origin.
\end{proof}


\section{Conjugacy classes of $\SU_{1,1}$}
\label{app:conj}

The Lie group $\SU_{1,1}$ is the set of matrices $S \in\SL_2\C$ satisfying
$S e_3 S^*=e_3$.  Two matrices $A,B\in \SU_{1,1}$ are
called {\it conjugate in\/} $\SL_2\C$ if there exists
a matrix $P\in \SL_2\C$ such that $B=P^{-1}AP$, and are called
{\it conjugate  in\/} $\SU_{1,1}$ if $B=P^{-1}AP$ for some
$P\in\SU_{1,1}$.

As in \eqref{eq:canonical-repr}  and \eqref{matrixR}, we set
 \begin{alignat*}{2}
  \Lambda_e(t)&:=
  \begin{pmatrix}
     e^{it} & 0 \\
      0     & e^{-it}
  \end{pmatrix}, \qquad
  &\Lambda_p(t) &:=
  \begin{pmatrix}
    1+i t & -it\\
    it & 1-it
  \end{pmatrix}, \\
  \Lambda_h(t) &:=
  \begin{pmatrix}
    \cosh t & \sinh t\\
    \sinh t & \cosh t
  \end{pmatrix}\quad
  &\quad \text{and}\quad
         R &:= \frac{1}{2}
  \begin{pmatrix}
   1 & \hphantom{-}1 \\ i & -i
  \end{pmatrix}
 \end{alignat*}
 for an arbitrary $t\in\R$.

\begin{theorem}\label{thm:canonical-su}
 A matrix $A\in\SU_{1,1}$ is conjugate in $\SU_{1,1}$
 to one of
 \begin{enumerate}
  \item $\Lambda_e(s)$ $(s\in (-\pi,\pi])$,
  \item $\pm\Lambda_p(t)$ or $\pm\Lambda_p(-t)$ $(t>0)$, or
  \item $\pm\Lambda_h(t)$ $(t>0)$.
 \end{enumerate}
\end{theorem}

\begin{remark}\label{rem:conj-neg}
 Though the matrices $\Lambda_e(s)$ and $\Lambda_e(-s)$ are
 conjugate in $\SL_2\C$, they are not conjugate in $\SU_{1,1}$
 if $s\not\equiv 0\pmod{2\pi}$.
 That is, for any elliptic matrix $A\in\SU_{1,1}$, there exists a
 {\em unique\/} real number $t\in (-\pi,\pi]$ such that
 $A$ and $\Lambda_e(s)$ are conjugate in $\SU_{1,1}$.
\end{remark}
\begin{remark}\label{rem:parab-conj}
 On the other hand,  $\Lambda_p(t_1)$ and $\Lambda_p(t_2)$ 
 ($t_1,t_2\neq 0$)  are conjugate in  $\SU_{1,1}$
 if and only if $t_1t_2>0$.
 In fact, if $t_1\neq t_2$,
 $P\Lambda_p(t_1)  P^{-1} = \Lambda_p(t_2)$ ($P\in\SU_{1,1}$) holds 
 if and only if 
 \[
     P = \pm \begin{pmatrix}
          a & \bar b \\
          b & \bar a
         \end{pmatrix},\qquad
         a = \cosh s + i u, ~ b=\sinh s + i u,
 \]
 where $s = \log\sqrt{t_2/t_1}\in\R$ and $u\in\R$.

 In particular, the sign of $t$ in $\Lambda_p(t)$
 is invariant under such a conjugation.
 Though $\Lambda_p(t)$ ($t\in\R\setminus\{0\}$)
 is conjugate with $\Lambda_p(1)$ or $\Lambda_p(-1)$,
 we choose various values of $t$ in this paper for the sake
 of convenience.
\end{remark}

\begin{remark}\label{rem:hyp-conj}
 Since
 \[
     \begin{pmatrix}
        i & \hphantom{-}0 \\
        0 & -i
     \end{pmatrix}
     \begin{pmatrix}
        \cosh t & \sinh t \\
        \sinh t & \cosh t
     \end{pmatrix}
     \begin{pmatrix}
        -i & 0 \\
        \hphantom{-}0 & i
     \end{pmatrix}=
     \begin{pmatrix}
        \cosh t & -\sinh t \\
        -\sinh t & \cosh t
     \end{pmatrix}
 \]
 $\Lambda_h(t)$ and $\Lambda_h(-t)$ are conjugate in $\SU_{1,1}$.
\end{remark}
\vspace{\baselineskip}
To show Theorem~\ref{thm:canonical-su},
we use the following group isomorphism:
\[
  \rho\colon{}\SL_2\R\ni X \longmapsto R^{-1}XR \in \SU_{1,1}.
\]
Note that $\Lambda_e(t), \Lambda_p(t)$, and $\Lambda_h(t)$
are the images of
\[
  \begin{pmatrix}
  \hphantom{-}\cos t & \sin t \\ -\sin t & \cos t \end{pmatrix},
  \quad
  \begin{pmatrix} 1 & 2t \\ 0  & 1 \end{pmatrix},\quad
  \text{and}\quad
  \begin{pmatrix} e^t & 0 \\ 0  & e^{-t} \end{pmatrix}
\]
 respectively, by $\rho$.

\begin{lemma}\label{key-L}
 Let $A$ and $B$ be $2\times 2$ \emph{real} matrices that are
 conjugate in $\SL_2\C$.
 Then $A$ and either $B$ or $e_3Be_3$ are conjugate
 in $\SL_2\R$.
\end{lemma}
\begin{proof}
 By assumption, there is a $\widetilde P\in \SL_2\C$
 with $A\widetilde P=\widetilde PB$.
 We set $\widetilde P=U+iV$, for real matrices $U$, $V$.
 Then $AU=UB$, $AV=VB$ and
 \[
    A(U+tV)=(U+tV)B \qquad\text{for any } t\in \R.
 \]
 If $\det(U+tV)$ vanishes identically for all $t \in \R$,
 holomorphicity of $\C \ni t \mapsto \det(U+tV)$ yields
 that $\det(U+iV) = \det \widetilde P=0$, a contradiction.
 Thus for some $t_0\in \R$, $\det(U+t_0V)\ne 0$, and
 then $(U+t_0V)^{-1}A(U+t_0V)=B$.
 If $\det (U+t_0V)>0$, we set
 $P=(U+t_0V)/\sqrt{\det(U+t_0V)}\in\SL_2\R$,
 giving $P^{-1}AP=B$.  If $\det(U+t_0V)<0$,
 we set $P=(U+t_0V)e_3/\sqrt{|\det(U+t_0V)|}\in\SL_2\R$,
 giving $P^{-1}AP=e_3Be_3$.
\end{proof}
\begin{proof}[Proof of Theorem~\ref{thm:canonical-su}]
 Let $A\in\SU_{1,1}$ and $\widetilde A:=\rho^{-1}(A)\in\SL_2\R$.

 If the eigenvalues of $A$ are not real numbers,
 they are written as $\{e^{it},e^{-it}\}$, where
 $t\in (-\pi,0) \cup (0,\pi)$.
 In this case, $\widetilde A$ is conjugate in $\SL_2 \C $ to
 $B_e := R \Lambda_e(t) R^{-1}$.  Hence by Lemma~\ref{key-L},
 $\widetilde A$ is conjugate in $\SL_2\R$ to $B_e$ or $e_3B_e e_3$.
 Thus, $A = \rho(\widetilde A)$ is conjugate in $\SU_{1,1}$ to
 $\rho(B_e)= \Lambda_e(t)$ or $\rho(e_3 B_e e_3)= \Lambda_e(-t)$.

 If the eigenvalues of $A$ are $\{\varepsilon,\varepsilon\}$
 ($\varepsilon=\{-1,1\}$) and $A\neq \varepsilon\id$,
 $\widetilde A$ is conjugate in $\SL_2 \C$ to
 $B_p := \varepsilon R \Lambda_p(t) R^{-1}$ for any $t\in\R_+$.
 Hence $\widetilde A$ is conjugate in $\SL_2\R$
 to either $B_p$ or $e_3 B_p e_3$.
 Thus, $A$ is conjugate in $\SU_{1,1}$ to
 $\rho(B_p)= \varepsilon \Lambda_p(t)$ or 
 $\rho(e_3 B_p e_3)=\varepsilon \Lambda_p(-t)$.
 As mentioned in Remark~\ref{rem:parab-conj},
 $\Lambda_p(u)$ for $u \in \R \setminus \{ 0 \}$
 is conjugate in $\SU_{1,1}$ to
 $\Lambda_p(\sign u)=\Lambda_p(\varepsilon)$.

 If the eigenvalues of $A$ 
 are two distinct real numbers, they are represented
 as $\{\varepsilon e^t,\varepsilon e^{-t}\}$, 
 where $t\in\R_+$ and $\varepsilon\in\{-1,1\}$.
 Thus, $\widetilde A$ is conjugate in $\SL_2\C$ to the diagonal matrix
 $B_h:=\epsilon R \Lambda_h(t)R^{-1}$.
 Hence, similarly to the first case, $A$ is conjugate in
 $\SU_{1,1}$ to $\rho(B_h)= \varepsilon \Lambda_h(t)$.
\end{proof}

\end{document}